\documentclass[ envcountsect,]{svjour3}
\smartqed
\usepackage{mystyle}
\usepackage{siunitx}
\newcommand{\beps}{\bvarepsilon}
\newcommand{\bepsp}{\bvarepsilon^\text{p}}

\newcommand{\dOm}{\,\text{d}\Omega}
\newcommand{\dev}{\operatorname{dev}}
\newcommand{\bCel}{\boldsymbol{C}_\text{el}}
\newcommand{\dz}{\,\text{dz}}
\newcommand{\Gc}{G_\text{c}}
\newcommand{\Gceff}{G_\text{c}^\text{eff}}
\graphicspath{{./pic/}}
\usepackage{breakcites}

\begin{document}

\title{Applications of conic programming in non-smooth mechanics}

\author{Jeremy Bleyer}

\institute{Jeremy Bleyer \at
             Laboratoire Navier, ENPC, Univ Gustave Eiffel, CNRS \\
              6-8 av. Blaise Pascal, Cité Descartes\\
              77455 Marne-la-vallée, France\\
              jeremy.bleyer@enpc.fr
}

\date{Received: date / Accepted: date}

\maketitle

\begin{abstract}
In the field of nonlinear mechanics, many challenging problems (e.g. plasticity, contact, masonry structures, nonlinear membranes) turn out to be expressible as conic programs. In general, such problems are non-smooth in nature (plasticity condition, unilateral condition, etc.), which makes their numerical resolution through standard Newton methods quite difficult. Their formulation as conic programs alleviates this difficulty since large-scale conic optimization problems can now be solved in a very robust and efficient manner, thanks to the development of dedicated interior-point algorithms. In this contribution, we review old and novel formulations of various non-smooth mechanics problems including associated plasticity with nonlinear hardening, nonlinear membranes, minimal crack surfaces and visco-plastic fluid flows.
\end{abstract}


\section{Introduction}
In full generality, a conic program denotes the problem of optimizing a linear function over the intersection of an affine subspace and a convex cone $\Kk$, for instance:
\begin{equation}
\begin{array}{rl} \displaystyle{\min_{\bx}} & \bc\T\bx \\
\text{s.t.} & \bA\bx = \bb \\
& \bx\in \Kk\end{array} \label{eq:conic-program}
\end{equation}

It turns out that every convex optimization problem consisting of minimizing a convex function $f(\bx)$ over inequality constraints defined using convex functions $f_i(\bx)$ can be expressed as a conic program of the form \eqref{eq:conic-program}.

However in practice, efficient solving algorithms are available for conic programs associated with a specific class of cones only, the so-called \textit{magic family} \cite{juditsky2021well}. One can mention in particular:
\begin{itemize}
\item positive orthants $\RR_+^m$;
\item Lorentz quadratic (or second-order) cones:
\begin{equation}
\Qq_{m}=\{(\bx,t)\in\RR^{m-1}\times\RR_+ \text{ s.t. } \|\bx\|_2 \leq t\}
\end{equation}
\item semi-definite cones $\Ss_{m}^+$, the cone of semi-definite positive $m\times m$ symmetric matrices;
\item power cones $\Pp^\alpha_{m}$ parameterized by $\alpha$ s.t. $0<\alpha < 1$:
\begin{equation}
\Pp^\alpha_m = \{\mathbf{z} \in \RR^m \text{ s.t. } \mathbf{z}=(z_0, z_1, \bar{\mathbf{z}}) \text{ and } z_0^\alpha z_1^{1-\alpha} \geq \|\bar{\mathbf{z}}\|_2, \:\: z_0,z_1\geq 0\}
\end{equation}
\item exponential cones:
\begin{equation}
\Kk_\text{exp} = \{\mathbf{z} \in \RR^{3} \text{ s.t. } \mathbf{z}=(z_0, z_1, z_2) \text{ and } z_0 \geq z_1\exp(z_2/z_1), \:\: z_0,z_1\geq 0\}
\end{equation}
\end{itemize}

Conic programs containing only positive orthants yield linear programs (LP), those containing quadratic Lorentz cones yield second-order cone programs (SOCP) and programs containing semi-definite cones yield semi-definite programs (SDP). An important property shared by these cones is that they are self-dual which is a corner-stone for the development of efficient interior-point (IP) solvers \cite{andersen2003implementing}. Non-self dual power and exponential cones \cite{dahl2021primal} were recently introduced in solvers such as Mosek for instance \cite{mosek}, increasing even more the modeling capabilities using such simple elementary bricks. For practical reasons, we will therefore refer to conic programs as those using cones belonging to these categories and which are efficiently solvable in practice, rather than the general abstract setting.

As a result, since every convex program can be turned into a conic format, we will also distinguish convex functions which are representable using such cones from generic abstract convex functions. Despite limiting ourselves to this family of cones, it turns out that it is sufficient to express a wide range of convex optimization functions appearing in practice \cite{mosekcookbook}.\\

Convex variational problems frequently arise in the modelling of many physical systems, in particular in the field of solid mechanics. As detailed later, the response of a material or a structure can often be described via an energy potential which encodes both balance and constitutive equations in a single variational principle. In a wide range of situations, such a potential happens to be convex. For instance, for linear elastic materials in small-strain, the latter is a simple positive-definite quadratic function. The associated optimality conditions therefore result in a simple linear variational equation. For more complex cases, the corresponding variational equation becomes non-linear which requires the use of a specific non-linear resolution method. Newton-type methods are most widely used and are very efficient when the system potential is smooth. However, in many important cases of applications, the underlying physics is described through some kind of threshold condition which correspond to non-smooth potentials. The seminal example is that of the obstacle problem where an elastic membrane deforms under some rigid obstacle position. The solution is therefore characterized by the coexistence of a contact and a no-contact zone in the geometrical domain. For this specific case, the associated problems turns out to be a bound-constrained quadratic problem. However, this system is still quite simple and active-set strategies are for instance very efficient to solve it. In more complex cases, such as for instance frictional contact \cite{wriggers2006computational}, elastoplasticity \cite{simo2006computational}, viscoplasticity \cite{dean2007numerical}, no-tension materials \cite{del1998limit}, the corresponding non-smooth variational problem is much more difficult to solve and tailored numerical methods must be developed (fixed point iterations, semi-smooth Newton methods, Augmented Lagrangian algorithms, etc.). 

Another approach is to consider the corresponding optimization problem as a conic programming problem and solve it using IP solvers, provided that the corresponding potential can be represented using the "magic cone" family. This approach has been mainly undertaken to solve limit analysis problems which are inherently formulated as constrained optimization problems which look for the maximal load factor a structure can sustain under equilibrium and convex yield conditions. Although alternative approaches have been proposed, the conic programming approach has emerged as the method of choice for solving large-scale limit analysis problems with applications ranging from soil mechanics \cite{krabbenhoft2008three} to steel construction \cite{el2020elastoplastic}, reinforced concrete \cite{vincent2018yield} or masonry structures \cite{portioli2014limit}. Following the success of such methods in a limit analysis setting, some contributions explored their application to elastoplastic problems \cite{krabbenhoft2007interior,krabbenhoft2007formulation}, including a recent extension towards non-convex finite strain plasticity \cite{el2021extending}. Topology optimization and plastic design of structures \cite{strang1983hencky} have also been formulated either as LP programs for trusses \cite{gilbert2003layout} or generic conic programs for solids \cite{mourad2021topology}. In the field of non-Newtonian fluids, viscoplastic (or yield stress) present the peculiarity of flowing like a liquid only when the stress reaches some yield stress limit. Below this yield stress, it behaves as a rigid solid. Interior-point algorithms based on conic programming have also been proposed as an efficient alternative to Augmented Lagrangian approaches which were commonly used in this community \cite{bleyer2015efficient,bleyer2018advances}. Finally, let us finish by mentioning the reference book on convex optimization in non-smooth mechanics \cite{kanno2011nonsmooth} which, in addition to considering frictional contact, plasticity and no-tension masonry structures, also deals with robust compliance optimization, cable networks, nonlinear membranes, etc.\\

Finally, we can see that conic programming formulations arise for many different non-smooth mechanical problems and become attractive due to the use of efficient IP solvers for their resolution. It is therefore not possible to give a complete overview of all these applications. Rather, we decide to select in this manuscript a few topics, some being old like elastoplasticity, others being original, and give new insights on modern conic programming formulations of such non-smooth mechanics problems.

The manuscript is organized as follows. Section \ref{sec:conic-repr} discusses generic conic representation of convex functions and variational problems and presents the numerical framework used to produce the simulation results. Section \ref{sec:evolution-problems} presents a generic incremental variational formulation for evolution problems of standard generalized materials. This setting encompasses a wide range of material behaviour, either rate-dependent or rate-independent. Section \ref{sec:plasticity} then applies the previous incremental principle to a hardening plasticity example. In particular, we consider an exponential hardening model which is formulated using recently developed exponential cones in IP solvers. Section \ref{sec:minimal-crack} deals with the determination of minimal crack surfaces in heterogeneous materials. It draws some links with min-cut and max-flow problems. Both primal and dual problems are, for the first time, solved using a conic programming formulation. Section \ref{sec:membranes} then deals with nonlinear membranes as in \cite{kanno2011nonsmooth} but considers a more generic formulation based on Ogden-type potentials which are reformulated using the newly introduced power cones in IP solvers. Finally, we finish with viscoplastic fluids flowing into the shallow geometry between two surfaces. We show that the dimensional reduction operation can be done from the underlying 3D fluid potential by retaining the conic-representation nature of the obtained potential. Again, this formulation is completely original to our knowledge. Section \ref{sec:conclusions} closes the manuscript with some conclusions and perspectives.

\begin{paragraph}{Notations}
In the following, vector and second-rank tensor will be represented using bold-face letters. We will use $\|\bx\|$ to denote the corresponding Euclidean norm of a quantity $\bx$ i.e. $\|\bx\|=\sqrt{\sum_i x_i^2}$ for a vector and $\|\bx\| = \sqrt{\sum_i\sum_j x_{ij}^2}$ for a second-rank tensor. The double product ":" will denote the inner product on second-rank tensors i.e $\bx:\by = \sum_i\sum_j x_{ij}y_{ji}$. Fourth-rank tensors will be represented with blackboard capital letters. 
\end{paragraph}

\section{Conic representation of convex functions and variational problems}\label{sec:conic-repr}
 In this work, given a convex function $f(\bx)$, we say that it is conic-representable if it can be written in the following form:
\begin{equation}
\begin{array}{rl}
\displaystyle{f(\bx) = \inf_{\by}} & \bc_x\T\bx+\bc_y\T\by \\
\text{s.t.} & \bb_l \leq \bA\bx + \bB \by\leq \bb_u \\
& \by \in \Kk
\end{array}\label{eq:conic-representation}
\end{equation}
where $\bc_x$, $\bc_y$, $\bb_l$, $\bb_u$ are given vectors of appropriate size, $\bA$ and $\bB$ are given matrix and $\Kk = \Kk_1\times\ldots \times \Kk_n$ is a product of cones where each $\Kk_i$ belongs to one of the previously mentioned class of cones. Note that the above representation is not necessarily unique. However, as soon as one equivalent form of \eqref{eq:conic-representation} is available, minimizing $f(\bx)$ can be done using dedicated IP solvers. As mentioned in \cite{juditsky2021well}, any given convex function does not necessarily possess an obvious conic reformulation and specific calculus rules must therefore be used in order to obtain such a representation. This process can be automated using Disciplined Convex Programming, as done for instance in the \texttt{CVX} software \cite{grant2014cvx}. Without being exhaustive, we give in the following several examples of simple conic-representable functions and some convexity-preserving operations which also maintain conic representability from that of the original function(s).
 
\subsection{Simple conic-representable functions}
Many interesting conic-representable functions can be obtained using simple building blocks such as:
\begin{itemize}
\item affine functions
\item positive definite quadratic functions
\item $p$-norms with $p\geq 1$: $p=2$ requires quadratic cones, $p=1$ or $\infty$ requires linear inequalities and $p\notin \{1,2,\infty\}$ can be expressed using generic power cones.
\item logarithm and exponential functions using exponential cones
\item maximum and minimum principal values using SDP cones
\item and many more.
\end{itemize}
We refer for instance to \cite{mosekcookbook} for more details on conic formulations of other usual functions.

\subsection{Operations on functions}
Based on the above simple functions, more complex functions can be obtained using simple operations on functions which all preserve the conic-representability. For instance, the sum of two conic-representable $f_1(\bx)+f_2(\bx)$ functions is obviously also conic-representable. Without being exhaustive, other examples include:
\begin{itemize}
\item linear precomposition by a linear operator $A$:
\begin{equation}
(f\circ A)(\bx) = f(\bA\bx)
\end{equation}
\item epigraph:
\begin{equation}
(\bx, t) \in \operatorname{epi}_f \Leftrightarrow f(\bx) \leq t
\end{equation}
\item perspective:
\begin{equation}
\operatorname{persp}_f(\bx, t) = tf(\bx/t)
\end{equation}
\item inf-convolution $f_1 \square f_2$ between two functions $f_1$,$f_2$:
\begin{equation}
\begin{array}{rl}
\displaystyle{(f_1 \square f_2) (\bx) = \inf_{\bx_1,\bx_2}} & f_1(\bx_1)+f_2(\bx_2)\\
\text{s.t.} & \bx = \bx_1+\bx_2
\end{array}
\end{equation}
\item the marginal $f\setminus \bA$ of $f$ through a linear operator $\bA$ \cite{fitzpatrick2001conjugates}:
\begin{equation}
\begin{array}{rl}
\displaystyle{(f \setminus \bA) (\bx) = \inf_{\by}} & f(\by)\\
\text{s.t.} & \bx = \bA\by
\end{array}
\end{equation}
\item a generalized marginal operator through a set of $n$ linear operators $\bA_i$ and positive weights $c_i \geq 0$:
\begin{equation}
\begin{array}{rl}
\displaystyle{(f \setminus \{\bA_i,c_i\}) (\bx) = \inf_{\by_1,\ldots,\by_n}} & \displaystyle{\sum_{i=1}^n c_i f(\by_i)}\\
\text{s.t.} & \displaystyle{\bx = \sum_{i=1}^n \bA_i\by_i}
\end{array} \label{eq:marginals}
\end{equation}
\item convex conjugate:
\begin{equation}
f^*(\by) = \sup_{\bx} \by\T\bx-f(\bx)
\end{equation}
\end{itemize}

\subsection{Conic-representable variational problems}

In mechanics, many problems benefit from a variational principle which consists in minimizing some potential $\Psi(\bu)$ with respect to some mechanical field $\bu$. Often this potential can be expressed from a corresponding potential density $\psi(\bu)$ per unit volume. The resulting variational problem formulated on a domain $\Omega$ generally reads as:
\begin{equation}
\min_{\bu\in V} \int_\Omega \psi(\bu)\dOm  \label{eq:var-princ-cont}
\end{equation}
where $V$ is some affine, or at least convex, space gathering various constraints on $\bu$ e.g. boundary conditions, physical bounds, etc.\\

When such a potential density $\psi$ turns out to be convex, the possibility of using conic programming to solve the problem therefore amounts to knowing whether $\psi$ is conic-representable or not using the conic magic family. As regards numerical resolution of the continuous variational problem \eqref{eq:var-princ-cont}, it is usually transformed into a discrete finite-dimensional variational problem by considering a finite-dimensional subspace $V_h\subset V$ and by evaluating the corresponding integral through a numerical quadrature such as:
\begin{equation}
\min_{\bar\bu=(\bu_i)\in V_h} \sum_{g=1}^G \omega_i \psi(\bL_g\bar\bu)  \label{eq:var-princ-discr}
\end{equation}
where $\omega_g$ are positive quadrature weights and $\bL_g$ are linear operators relating the vector of unknowns $\bar\bu=(\bu_i)\in V_h$ to the value of the continuous field $\bu$ at the corresponding quadrature point.

\subsection{Numerical tools}
In the present work, we heavily rely on the FEniCS software package \cite{logg2012automated} for discretizing variational problems using the finite-element method. FEniCS provides various types of interpolation ranging from standard Lagrange polynomials, either in a continuous Galerkin or a discontinuous Galerkin form but also Raviart-Thomas elements, Brezzi-Douglas-Marini elements, etc. Concerning the resolution of discrete conic programs, we will make use of Mosek v9.0. The formulation of conic variational programs within the FEniCS environment is performed using the \texttt{fenics\_optim} software package which is freely available \cite{fenics_optim_2020}. \texttt{fenics\_optim} is a domain-specific language which enables easy definition of convex variational problems via conic representation of functions/constraints and combination of such using pre-implemented operators (e.g. perspective, inf-convolution, epigraph, etc.). For more details on the interaction between conic programming formulations and finite-element discretization, we refer the reader to the paper \cite{bleyer2020automating} which complements the package documentation.

Since spatial discretization is not the primary focus of the present work, the following variational problems will all be written in a continuous setting such as \eqref{eq:var-princ-cont}. It is however implied that their concrete implementation and resolution is performed using a specific discretization such as \eqref{eq:var-princ-discr} through the choice of an appropriate finite-element function space $V_h$ and a quadrature rule.

\section{Conic programming formulation of evolution problems for generalized standard materials}\label{sec:evolution-problems}

\subsection{Generalized standard materials}
We consider here the context of generalized standard materials \cite{halphen1975materiaux} which are characterized by a set of state variables, a free energy potential and a dissipation pseudo-potential. This thermodynamic approach benefits from a variational formulation which satisfies the first and second fundamental thermodynamic principles under some appropriate convexity assumptions for the free energy and dissipation pseudo-potential. 

Ignoring temperature effects for simplicity and restricting to a small strain setting, we assume that state variables consist here of the total linearized strain $\beps$ and a set of internal state variables $\balpha$ describing the irreversible behavior of the material. The material free energy density is denoted by $\psi(\beps,\balpha)$. We assume here that $\psi$ is convex. Finally, we introduce $\phi(\dot{\beps},\dot{\balpha})$ the pseudo-potential which is also assumed to be a convex function of $(\dot{\beps},\dot{\balpha})$. In general, one could assume that $\phi$ also depends on $(\beps,\balpha)$. For the sake of simplicity, we assume here that this is not the case.

\subsection{Incremental pseudo-potential}
Consider now a time increment $[t_n,t_{n+1}]$ for which the state $(\beps_n,\balpha_n)$ at time $t_n$ is known. 

It can be shown that the solution at time $t_{n+1}$ in terms of displacement $\bu$ and state variables can be obtained as the solution to the following minimization principle \cite{ortiz1999variational,mielke2005evolution}:
\begin{equation}
\begin{array}{rll} \displaystyle{(\bu_{n+1},\beps_{n+1},\balpha_{n+1}) =\argmin_{
(\bu,\beps,\balpha) }} &
 \displaystyle{\int_{t_n}^{t_{n+1}}\int_{\Omega}(\dot{\psi}(\beps,\balpha)+\phi(\dot{\beps},\dot{\balpha}))\dOm\,\text{dt}} \\
&\displaystyle{ - \int_{t_n}^{t_{n+1}} \mathcal{P}_\text{ext}(\dot{\bu})\,\text{dt}}\\
\text{s.t.} & (\bu,\beps)\in \text{KA}_{n+1} \end{array}
\end{equation}
where:
\begin{itemize}
\item $\mathcal{P}_\text{ext}(\dot{\bu})$ is the power of external loads which we assume to consist only of fixed body forces $\bf_{n+1}$ on the time step, so that:
\begin{equation}
\int_{t_n}^{t_{n+1}} \mathcal{P}_\text{ext}\,\text{dt} = \int_{t_n}^{t_{n+1}} \int_\Omega \bf_{n+1}\cdot\dot{\bu} \dOm=\int_\Omega \bf_{n+1}\cdot(\bu-\bu_n) \dOm
\end{equation}
\item $\text{KA}_{n+1}$ denotes the state of kinematically admissible fields defined as:
\begin{equation}
\text{KA}_{n+1} = \left\{(\bu,\beps) \text{ s.t. } \left|\begin{array}{ll}
\beps = \nabla^s \bu & \text{in }\Omega \\
\bu = \bar{\bu}_{n+1} & \text{on } \partial \Omega_\text{D} \end{array}\right.\right\}
\end{equation}
where $\bar{\bu}_{n+1}$ corresponds to imposed displacements at time $t_{n+1}$ on the Dirichlet boundary $\partial \Omega_\text{D}$.
\end{itemize} 

Further introducing a backward-Euler approximation for the state variables evolution : 
\begin{equation}
\dot{\beps}(t) \approx \dfrac{\beps-\beps_n}{\Delta t}\quad;\quad \dot{\balpha}(t) \approx \dfrac{\alpha-\alpha_n}{\Delta t}
\end{equation}
where $\Delta t =t_{n+1}-t_n$, we have:
\begin{equation}
\int_{t_n}^{t_{n+1}}\int_{\Omega}\phi(\dot{\beps},\dot{\balpha})\dOm\,\text{dt} \approx \int_{\Omega} \Delta t\phi\left(\dfrac{\beps-\beps_n}{\Delta t},\dfrac{\alpha-\alpha_n}{\Delta t}\right)\dOm 
\end{equation}

The above minimum principle can thus be replaced with:
\begin{equation}
\begin{array}{rll} \displaystyle{(\bu_{n+1},\beps_{n+1},\balpha_{n+1}) =\argmin_{
(\bu,\beps,\balpha) }} &
 \displaystyle{\int_{\Omega}\left(\psi(\beps,\balpha)-\psi(\beps_n,\balpha_n)\right)\dOm} \\
&\displaystyle{ + \int_{\Omega} \Delta t\phi\left(\dfrac{\beps-\beps_n}{\Delta t},\dfrac{\alpha-\alpha_n}{\Delta t}\right)\dOm} \\
&\displaystyle{-\int_\Omega \bf_{n+1}\cdot(\bu-\bu_n) \dOm}\\
\text{s.t.} & (\bu,\beps)\in \text{KA}_{n+1} \end{array} \label{eq:incr-principle-rate-dep}
\end{equation}

Ignoring the constant terms $\psi(\beps_n,\balpha_n)$ and $\bu_n$, \eqref{eq:incr-principle-rate-dep} becomes:
\begin{equation}
\begin{array}{rll} \displaystyle{(\bu_{n+1},\beps_{n+1},\balpha_{n+1}) =\argmin_{
(\bu,\beps,\balpha) }} &
 \displaystyle{\int_{\Omega}J(\beps,\balpha) \dOm - \int_\Omega \bf_{n+1}\cdot\bu \dOm}\\
\text{s.t.} & (\bu,\beps)\in \text{KA}_{n+1} \end{array} \label{eq:incr-principle}
\end{equation}
where:
\begin{itemize}
\item for rate-dependent materials:
\begin{equation}
J(\beps,\balpha) = \psi(\beps,\balpha)+\Delta t\phi\left(\dfrac{\beps-\beps_n}{\Delta t},\dfrac{\balpha-\balpha_n}{\Delta t}\right) 
\end{equation}
\item for rate-independent materials:
\begin{equation}
J(\beps,\balpha) = \psi(\beps,\balpha)+\phi\left(\beps-\beps_n,\balpha-\balpha_n\right) \quad \text{}
\end{equation}
since the pseudo-potential $\phi$ is, in addition, a positively homogeneous convex function in this specific case.
\end{itemize}
Finally, the incremental pseudo-potential $J$ is a convex function of $(\beps,\balpha)$ which makes problem \eqref{eq:incr-principle} a convex optimization problem. Provided that both $\psi$ and $\phi$ are conic representable, the resulting problem can be expressed as a conic program.

\section{Application to associated plasticity}\label{sec:plasticity}
\subsection{$J_2$-plasticity with nonlinear isotropic hardening}
We now particularize the previous generic formulation to the case of associated plasticity with nonlinear isotropic hardening. The internal state variables for this setting are the plastic strain $\bepsp$ and the cumulated equivalent plastic strain $p$. The free energy density $\psi$ consists of a stored elastic $\psi_\text{el}$ and hardening $\psi_\text{h}$ potentials:
\begin{equation}
\psi(\beps,\bepsp,p) = \psi_\text{el}(\beps-\bepsp)+\psi_\text{h}(p)
\end{equation}
We consider  isotropic linear elasticity:
\begin{equation}\psi_\text{el}(\beps^\text{el}) = \dfrac{1}{2}\beps^\text{el}:\CC:\beps^\text{el} = \dfrac{\kappa}{2}(\tr(\beps^\text{el}))^2 + \mu \dev(\beps^\text{el}):\dev(\beps^\text{el}) \label{eq:isotropic-elasticity}
\end{equation}
where $\kappa$ is the compressibility modulus, $\mu$ the shear modulus and $\dev(\beps^\text{el})$ the deviatoric elastic strain.

The hardening potential is assumed to be of exponential type as follows:
\begin{equation}
\psi_\text{h}(p) = (\sigma_u-\sigma_0)\left(p+\frac{1}{\omega}\exp(-\omega p)\right)
\end{equation}
where $\sigma_0$ (resp. $\sigma_u$) is the initial (resp. ultimate) yield strength and $\omega$ a saturation parameter. This potential defines the following hardening thermodynamic force:
\begin{equation}
R(p) = \dfrac{\partial \psi_\text{h}}{\partial p} =  (\sigma_u-\sigma_0)(1-\exp(-\omega p))
\end{equation}
which will increase the yield stress $R(p)$ from $\sigma_0$ at $p=0$ to $\sigma_u$ when $p\to\infty$.

Finally, we assume a $J_2$-plasticity dissipation pseudo-potential:
\begin{equation}
\phi(\dot{\bepsp}, \dot{p}) = \begin{cases}
\sqrt{\frac{2}{3}}\sigma_0\|\dot{\bepsp}\| & \text{if } \tr(\dot{\bepsp})=0 \\
+\infty & \text{otherwise}
\end{cases}\label{eq:J2-dissipation}
\end{equation}
which involves a plastic incompressibility constraint. However, as such, we are missing the link between the plastic strain and the cumulated equivalent plastic strain. Classically, one defines the equivalent plastic strain as follows:
\begin{equation}
p = \int_0^t \sqrt{\frac{2}{3}}\|\dot{\bepsp}\|\text{dt}
\end{equation}
or, equivalently in rate form:
\begin{equation}
\dot{p} = \sqrt{\frac{2}{3}}\|\dot{\bepsp}\| \label{eq:definition-dot-p}
\end{equation}

\subsection{Conic reformulations}

We now discuss the conic representation of the various convex functions  $\psi_\text{el},\psi_\text{h}$ and $\phi$ involved in the definition of $J(\beps,\bepsp,p)$.

Starting with the plastic dissipation \eqref{eq:J2-dissipation}, let us first mention that the constraint \eqref{eq:definition-dot-p} is non-convex. The proper constraint linking both internal state variable rates in a convex manner is to relax the equality and redefine $\phi$ using $\dot{p}$ as follows:
\begin{equation}
\phi(\dot{\bepsp}, \dot{p}) = \begin{cases}
\sigma_0 \dot{p} & \text{if } \tr(\dot{\bepsp})=0\:; \sqrt{\frac{2}{3}}\|\dot{\bepsp}\| \leq \dot{p} \\
+\infty & \text{otherwise}
\end{cases}\label{eq:J2-dissipation-conic}
\end{equation}
Formulation \eqref{eq:J2-dissipation-conic} corresponds to an epigraph form of \eqref{eq:J2-dissipation} which is readily expressed using a second-order cone constraint.

Considering now the elastic potential, the quadratic form \eqref{eq:isotropic-elasticity} must be expressed using a Cholesky factorization. Accounting for the fact that the spherical and deviatoric parts of a second-rank tensor define orthogonal subspaces, we can readily show that:
\begin{equation}
\psi_\text{el}(\beps^\text{el}) = \dfrac{1}{2}\|\QQ:\beps^\text{el}\|^2  \label{eq:isotropic-elasticity-factorized}
\end{equation}
where:
\begin{equation}
\QQ:\beps^\text{el} = \sqrt{3\kappa}\dfrac{1}{3}\tr(\beps^\text{el})\bI + \sqrt{2\mu} \dev(\beps^\text{el})
\end{equation}
which yields the following conic epigraph formulation:
\begin{equation}
\begin{array}{rl} \displaystyle{\psi_\text{el}(\beps^\text{el}) = \inf_{t, s, \by}} & t \\
\text{s.t.} & \by = \QQ:\beps^\text{el} \\
& \|\by\|^2 \leq 2t s \\
& s = 1
\end{array}  \label{eq:isotropic-elasticity-epigraph}
\end{equation}

Finally, the hardening potential term can be readily reformulated using an exponential cone as follows:
\begin{equation}
\exp(-\omega p) = \min \:r_0 \quad \text{s.t. } \exp(-\omega p) \leq r_0
\end{equation}

This non-linear constraint can be reformulated using an exponential cone $\Kk_\text{exp}$ as follows:
\begin{equation}
\exp(-\omega p) \leq r_0 \quad \Leftrightarrow \quad r_1=1, \:, r_2=-\omega p, \: (r_0, r_1, r_2)\in \Kk_\text{exp}
\end{equation}

\subsection{Numerical illustration}

\begin{figure}
\begin{center}
\includegraphics[width=0.8\textwidth]{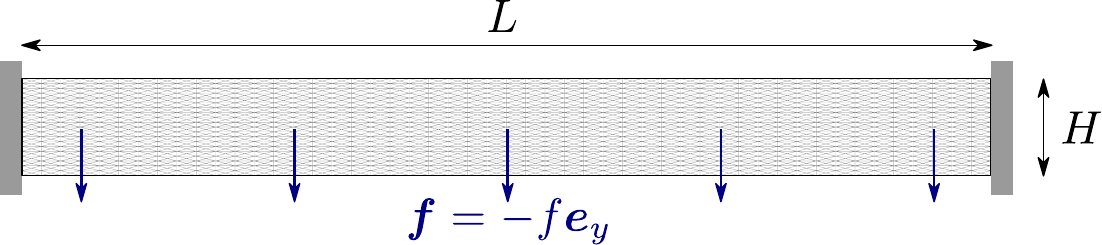}
\end{center}
\caption{A 2D beam clamped at both ends and subject to a vertical body force $f$.}\label{fig:beam-mesh}
\end{figure}

We consider a 2D rectangular domain of length $L=5$ and height $H=0.5$, fixed on both lateral extremities and subjected to a uniform downwards vertical body force $\bf=-f\be_y$ as illustrated in Figure \ref{fig:beam-mesh}. The domain is meshed with $50\times 20$ elements on the boundaries. The displacement field is discretized using a continuous quadratic Lagrange function space whereas internal state variables are represented using a discontinuous piecewise affine space. We consider a plane strain setting and use the following material properties $E=\SI{210}{GPa}$, $\nu=0.3$, $\sigma_0=\SI{450}{MPa}$, $\sigma_u=\SI{700}{MPa}$ and $\omega=50$.

The beam is first loaded by progressively increasing the body force from 0 to $f^+=\dfrac{2}{\sqrt{3}}\dfrac{4\sigma_u H}{L^2}$ which is the theoretical limit load using a beam theory model. Then, we perform a full unloading up to $f=0$. The loading stage is imposed using 10 load increments, the unloading stage being elastic using only one increment.

\begin{figure}
\begin{center}
\includegraphics[width=0.6\textwidth]{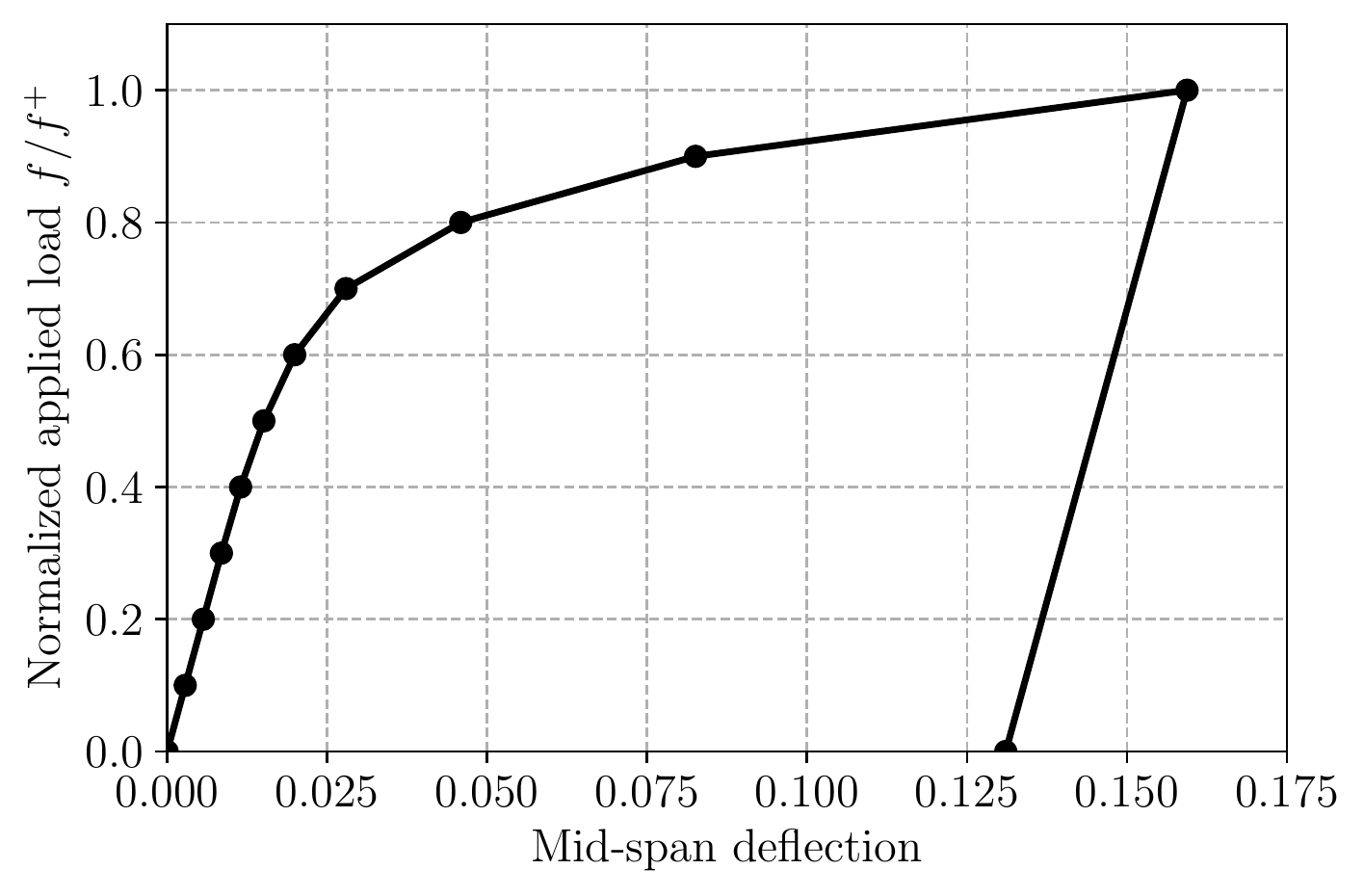}
\end{center}
\caption{Load vs. mid-span deflection evolution during plastic loading and elastic unloading.}\label{fig:elastoplastic-load-displ-10steps}
\end{figure}

\begin{figure}
\begin{center}
\includegraphics[width=0.8\textwidth]{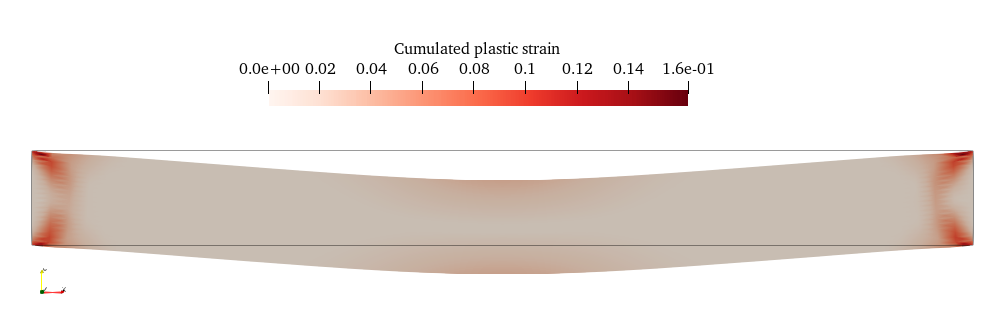}
\end{center}
\caption{Equivalent plastic strain distribution and deformed configuration at $f=f^+$.}\label{fig:elastoplastic-p-distribution}
\end{figure}

Figure \ref{fig:elastoplastic-load-displ-10steps} represents the evolution of the beam downwards vertical displacement at its mid-span center point $(L/2,H/2)$ as a function of the imposed loading. One can observe an initial elastic phase up to $f\approx 0.4f^+$ followed by a strongly non-linear hardening phase. Note that the load $f=f^+$ can still be supported by the structure due to the difference between a 2D model as here and a 1D beam theory solution. Note that, for the present 2D structure with the considered mesh, the ultimate load is found to be around $f\approx 1.1f^+$. As expected for plasticity problems, the unloading stage is indeed elastic and exhibits a permanent residual deflection. The distribution of equivalent plastic strain $p$ at $f=f^+$ has been represented in Figure \ref{fig:elastoplastic-p-distribution}. One can clearly observe the formation of plastic hinges near the clamped supports and a more diffuse plastic field at the beam mid-span. 

As regards numerical resolution statistics, let us point out that each instance of a single step elasto-plastic problem consists in  roughly 250,000 optimization variables, 170,000 linear constraints and 36,000 quadratic or exponential cones. Each resolution with Mosek v.9.0 took between 2.5 and 3.5 seconds (10 to 18 IP iterations) depending on the loading step (fully plastic steps near $f=f^+$ took a larger number of iterations than initial elastic steps). Overall, we find that the IP solver exhibits a very robust behaviour in terms of number of iterations with respect to the load step level or to the problem size (after mesh refinement for instance). 

\begin{figure}
\begin{center}
\includegraphics[width=0.8\textwidth]{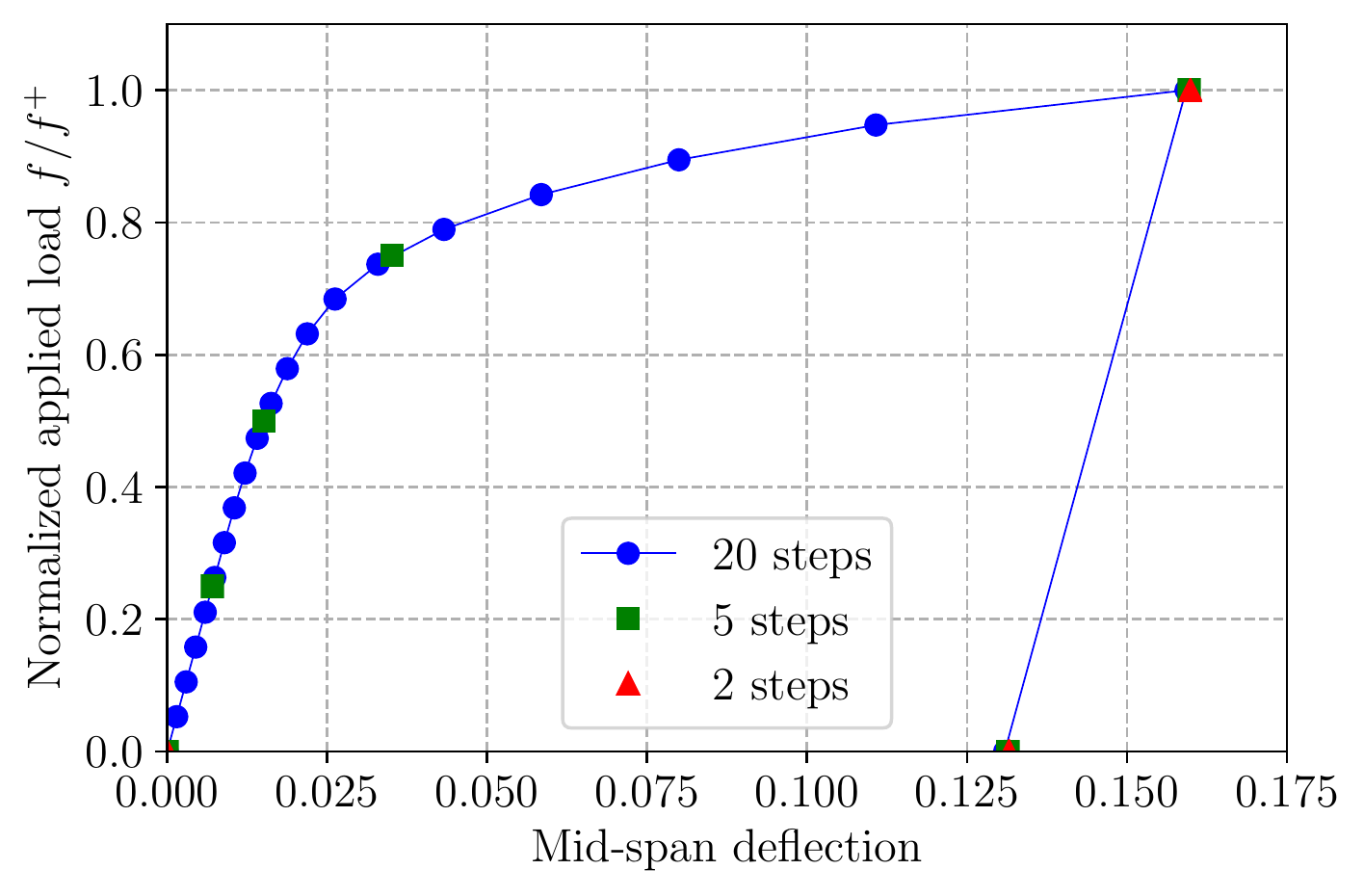}
\end{center}
\caption{Load-deflection evolution depending on the total number of load increments (only one used in the unloading phase).}\label{fig:elastoplasticity-increment-size}
\end{figure}

Although not being fully competitive against standard Newton methods in multi-step plasticity, the conic programming approach becomes interesting when much larger load steps are considered. For instance, Figure \ref{fig:elastoplasticity-increment-size} shows the same load-displacement curve as before using different numbers of load steps. Interestingly, the computed displacement at the ultimate state near collapse ($f=f^+$) is already very accurate using a single load step, see also \cite{krabbenhoft2007interior,el2020elastoplastic}. Again, solver robustness does not seem to be affected by the load step amplitude since the number of IP iterations remains very similar.

\section{Minimal crack surfaces}\label{sec:minimal-crack}
In this section, we consider the problem of computing the effective crack resistance of a heterogeneous medium with locally varying fracture energy $\Gc(\bx)$. We consider $\Omega$ to be some representative volume element (RVE) of the heterogeneous material. In \cite{braides1996homogenization}, a periodic homogenization result regarding the variational approach to fracture \cite{francfort1998revisiting} has been established. In particular, the effective fracture energy $\Gceff(\bn)$ associated with a crack of mean normal $\bn$ was explicitly characterized from the computation of minimal surfaces inside $\Omega$, weighted by the local fracture energy $\Gc(\bx)$. \cite{schneider2020fft} proposed a convex optimization formulation inspired by min-cut/max-flow problems in a periodic setting. More precisely, given a prescribed crack plane normal $\bn$, they consider the following variational problem:
\begin{equation}
\Gceff(\bn) = \inf_{\phi \in V} \dfrac{1}{|\Omega|}\int_\Omega \Gc(\bx)\|\nabla \phi + \bn\|_2 \dOm \label{eq:mincut-crack-surface-Gceff}
\end{equation}
where $V$ denotes the space of smooth scalar functions which are periodic over $\Omega$. This min-cut problem is known to have minimizers corresponding to $\Gc$-weighted periodic minimal surfaces.

\subsection{Illustrating example}
\begin{figure}
\begin{center}
\includegraphics[width=0.8\textwidth]{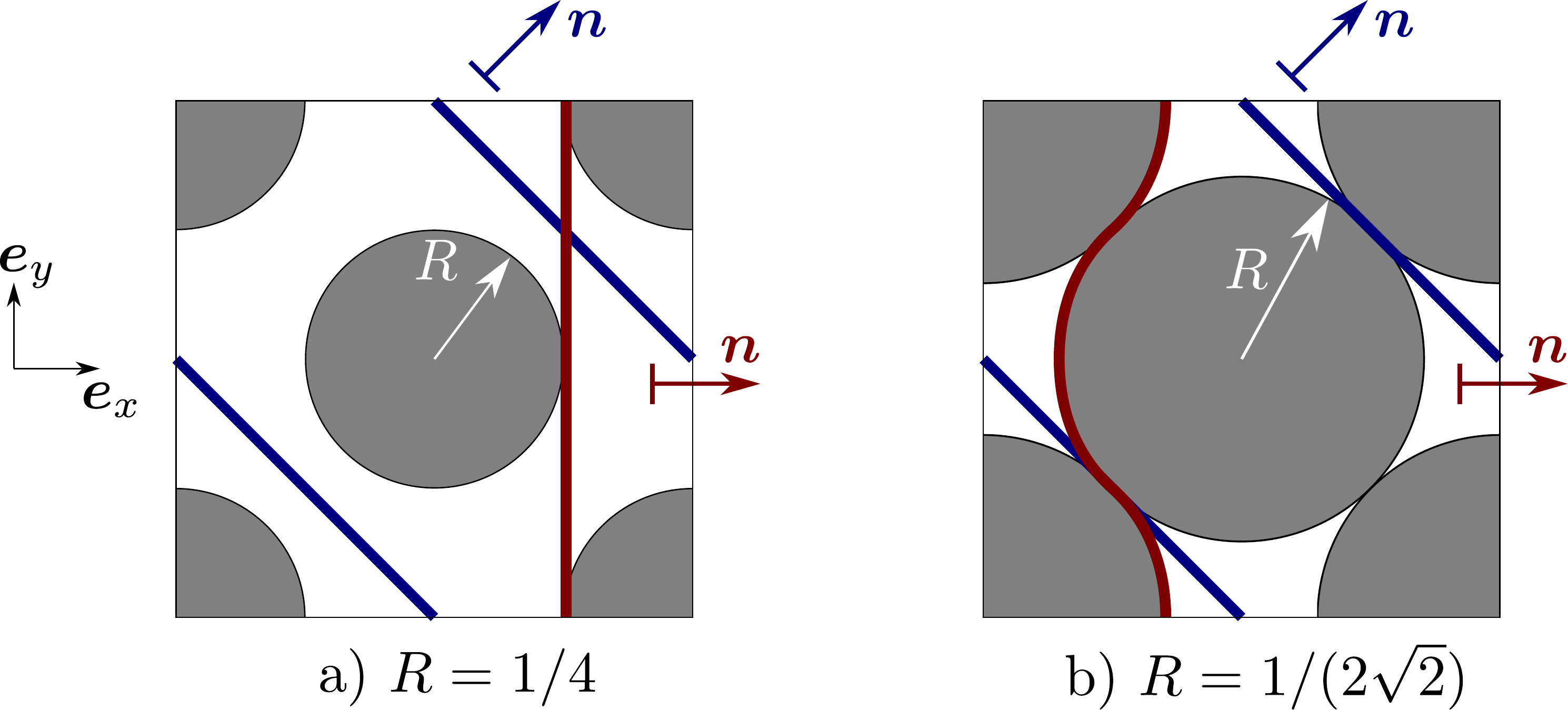}
\end{center}
\caption{Minimal crack surfaces for $0^\circ$ (red) and $45^\circ$ (blue) imposed crack plane orientation.}\label{fig:minimal-crack-surface-regular-inclusions}
\end{figure}

We consider a simple microstructure consisting of a periodic square unit cell and circular inclusions of radius $R \leq R_0=1/(2\sqrt{2})$ located at the cell center and at its four corners, see Figure \ref{fig:minimal-crack-surface-regular-inclusions}. $R_0$ denotes the maximum radius corresponding to each inclusion touching each other. The inclusion material possesses a fracture energy which is much larger than that of the matrix. Inclusions can thus be considered as infinitely resistant so that minimal crack surfaces will always pass through the matrix material only with fracture energy $\Gc$.

Let us first consider the effective fracture energy for cracks of normal $\be_x$ (or $\be_y$ for symmetry reasons). For small enough inclusions i.e. $R\leq 1/4$, there always exists a straight crack plane passing inside the matrix (Fig. \ref{fig:minimal-crack-surface-regular-inclusions}a). Thus, $\Gceff(\be_x)=\Gc$ for any $R\leq 1/4$. For larger $R$ up to $R=R_0$, the optimal path consists of a part of the circular inclusion border connected by a straight line between each inclusions. In the limit case $R=R_0$, the straight line vanishes and the total path length is $\pi R_0$ (Fig. \ref{fig:minimal-crack-surface-regular-inclusions}b). As a result, the effective fracture energy is $\Gceff(\be_x)=\dfrac{\pi}{2\sqrt{2}}\Gc \approx 1.11\Gc$.

Regarding crack planes oriented at $\pm 45^\circ$, for any $R$ up to $R_0$ there exists a $\pm 45^\circ$ degree line connecting the mid-point of the square domains as indicated in Figure \ref{fig:minimal-crack-surface-regular-inclusions}. As a result, the corresponding fracture energy will always be $\Gceff = \Gc$ for this orientation.

\subsection{Numerical computation}

The convex problem \eqref{eq:mincut-crack-surface-Gceff} is inherently non-smooth and hence difficult to solve in practice with standard Newton-methods. In \cite{schneider2020fft,ernesti2021fast}, first-order proximal algorithms have been proposed in conjunction with a Fast Fourier Transform-based discretization technique. In the present contribution, we can easily solve this problem using conic programming and a finite-element discretization. However, FFT-based proximal algorithms should be definitely more competitive for large scale 3D problems. 

Considering a triangular mesh and cell-wise constant values for $\Gc(\bx)$, we used a linear interpolation including periodic boundary conditions for the unknown field $\phi(\bx)$ in the min-cut problem \eqref{eq:mincut-crack-surface-Gceff}. The corresponding optimal value is an upper bound estimate to $\Gceff(\bn)$ which will converge upon mesh refinement.

Moreover, in order to assess the influence of spatial discretization on the computed fracture energy, we also solve the corresponding dual max-flow problem:
\begin{equation}
\begin{array}{rll}\displaystyle{\Gceff(\bn) = \sup_{\btau \in W}} & \displaystyle{\dfrac{1}{|\Omega|}\int_\Omega \btau\cdot\bn \dOm}  & \\
\text{s.t.} & \div \btau(\bx) = 0 & \forall \bx\in\Omega\\
& \|\btau(\bx)\|_2 \leq \Gc(\bx)  &
\end{array}
 \label{eq:maxflow-crack-surface-Gceff}
\end{equation}
where $W$ is a function space of smooth periodic vector fields with normal component continuity. Upon proper discretization of the above problem e.g. using a first-order Raviart-Thomas function space (see the Cheeger set example in \cite{bleyer2020automating}), the obtained optimal value can be shown to provide a lower bound to the exact effective crack energy $\Gceff(\bn)$. Solving both min-cut and max-flow problems with proper discretization will therefore furnish a bracketing of the exact value.

\subsection{Results on the illustrative example}
\begin{figure}
\begin{center}
\includegraphics[width=0.8\textwidth]{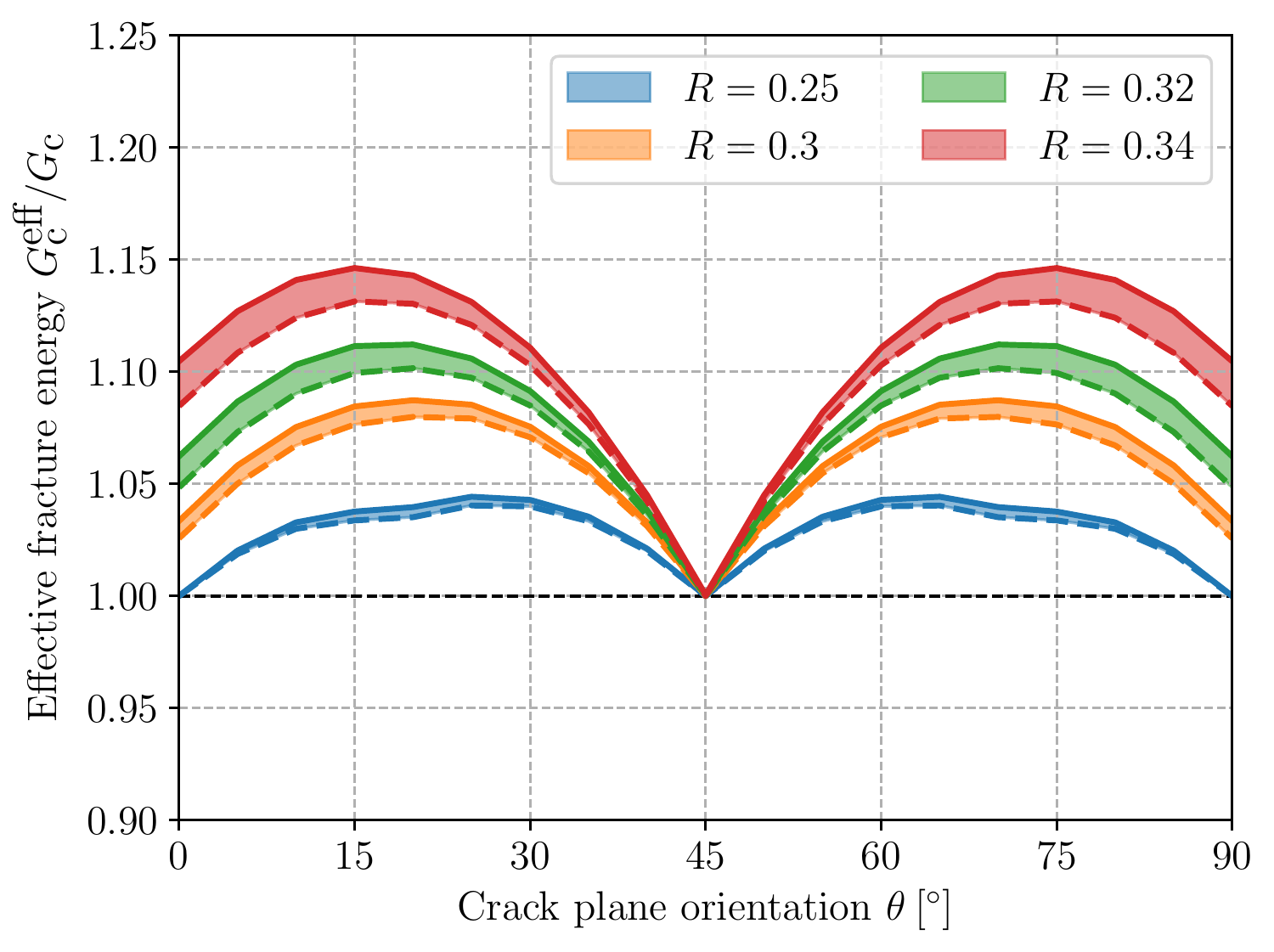}
\end{center}
\caption{Evolution of $\Gceff(\bn)$ for the microstructure in Fig. \ref{fig:minimal-crack-surface-regular-inclusions} as a function of the crack plane orientation $\bn(\theta)$. Colored regions are delimited by the min-cut upper bound (solid lines) and the max-flow lower bound (dashed lines).}\label{fig:minimal-crack-surface-regular-inclusions-results}
\end{figure}

\begin{figure}
\begin{center}
\begin{subfigure}{0.24\textwidth}
\includegraphics[width=\textwidth]{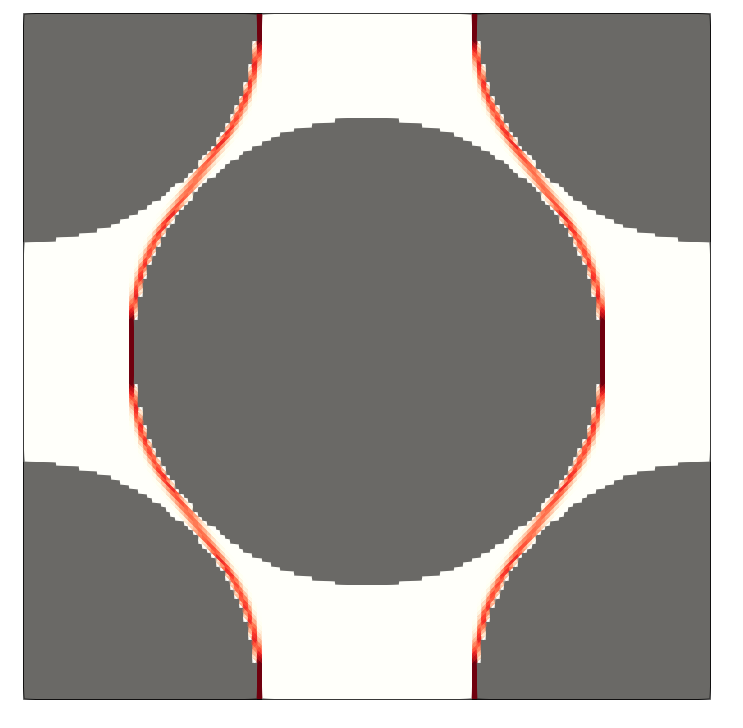}
\caption{$\theta=0^\circ$}
\end{subfigure}
\hfill
\begin{subfigure}{0.24\textwidth}
\includegraphics[width=\textwidth]{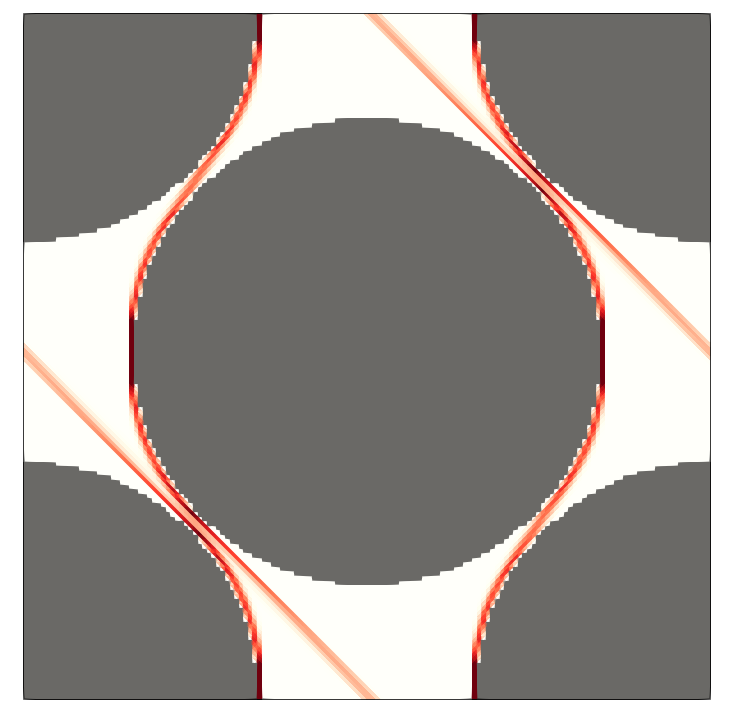}
\caption{$\theta=15^\circ$}
\end{subfigure}
\hfill
\begin{subfigure}{0.24\textwidth}
\includegraphics[width=\textwidth]{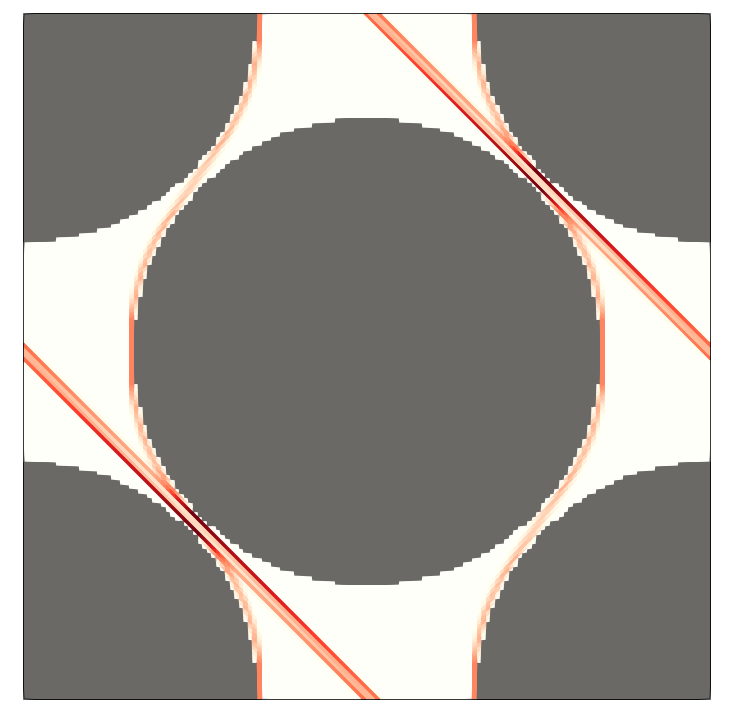}
\caption{$\theta=30^\circ$}
\end{subfigure}
\hfill
\begin{subfigure}{0.24\textwidth}
\includegraphics[width=\textwidth]{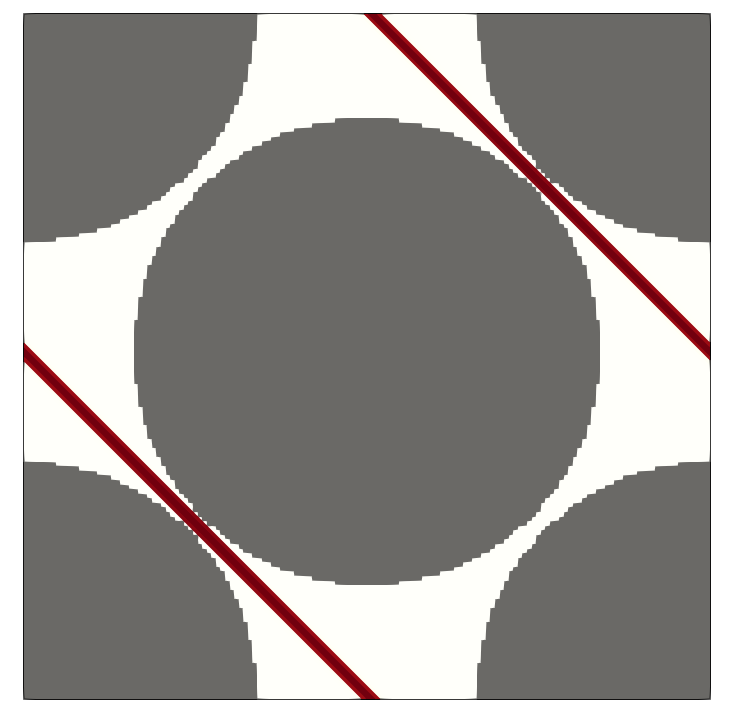}
\caption{$\theta=45^\circ$}
\end{subfigure}
\end{center}
\caption{Crack density field for $R=0.34$}\label{fig:minimal-crack-surface-regular-inclusions-crack-density}
\end{figure}

Figure \ref{fig:minimal-crack-surface-regular-inclusions-results} represents the evolution of $\Gceff(\bn)$ for the microstructure in Fig. \ref{fig:minimal-crack-surface-regular-inclusions} and $\bn=\cos\theta \be_x+\sin\theta \be_y$. Colored regions denote the bracketing obtained on a $150\times 150$ mesh when solving both min-cut and max-flow problems \eqref{eq:mincut-crack-surface-Gceff} and \eqref{eq:maxflow-crack-surface-Gceff}. First, we indeed obtain $\Gceff=\Gc$ for any $R\leq R_0$ at $\theta=\pm 45^\circ$ and at $\theta = 0^\circ$ or $90^\circ$ for $R\leq 1/4$. When $R$ approaches $R_0\approx 0.3535$, we have $\Gceff \approx \pi/R_0 \approx 1.11$ as expected. We also notice that the relative error between both bounds increases when $R$ approaches the limit $R_0$. The optimal crack density fields $\|\nabla\phi+\bn\|$ have been represented in Figure \ref{fig:minimal-crack-surface-regular-inclusions-crack-density}. Let us note that minimal crack surfaces are not necessarily unique. For instance, in Fig. \ref{fig:minimal-crack-surface-regular-inclusions}b, the mirror path on the right of the central inclusion is also a minimizer for the case $\theta=0^\circ$ as well as any combination of the left and right paths. The minimizer of the discrete problem seems to select a symmetric solution (Fig. \ref{fig:minimal-crack-surface-regular-inclusions-crack-density}a) very close to the expected solution of Fig. \ref{fig:minimal-crack-surface-regular-inclusions}b. For $\theta=45^\circ$, we also obtain the expected solution. Minimizers for intermediate angles $\theta$ seem to be obtained as a weighted combination of the $0^\circ$ and $45^\circ$ paths, which results in an increased effective fracture energy.

\subsection{Application to a microstructure with random inclusions}

\begin{figure}
\begin{center}
\begin{subfigure}{0.44\textwidth}
\includegraphics[width=\textwidth]{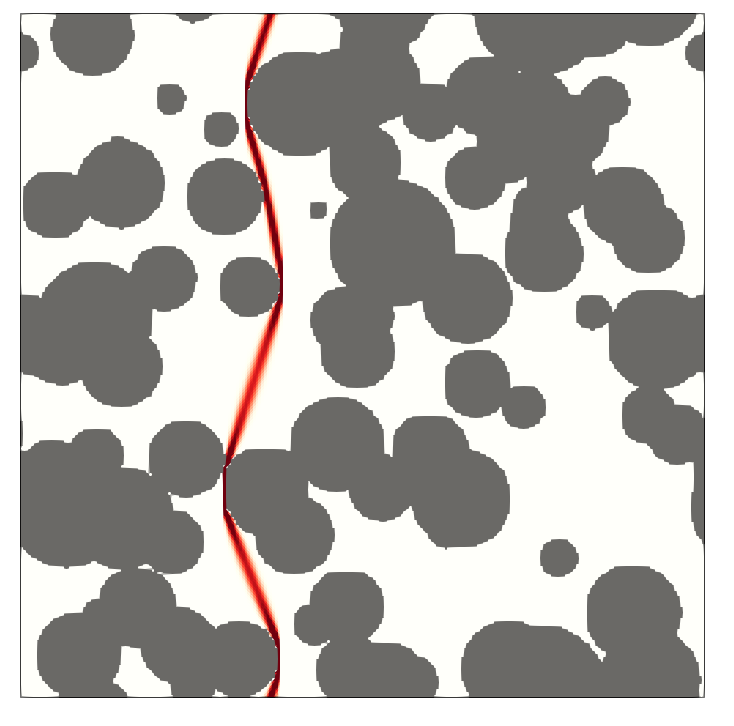}
\caption{$\theta=0^\circ$, $\Gceff=1.048\Gc$}
\end{subfigure}
\hfill
\begin{subfigure}{0.44\textwidth}
\includegraphics[width=\textwidth]{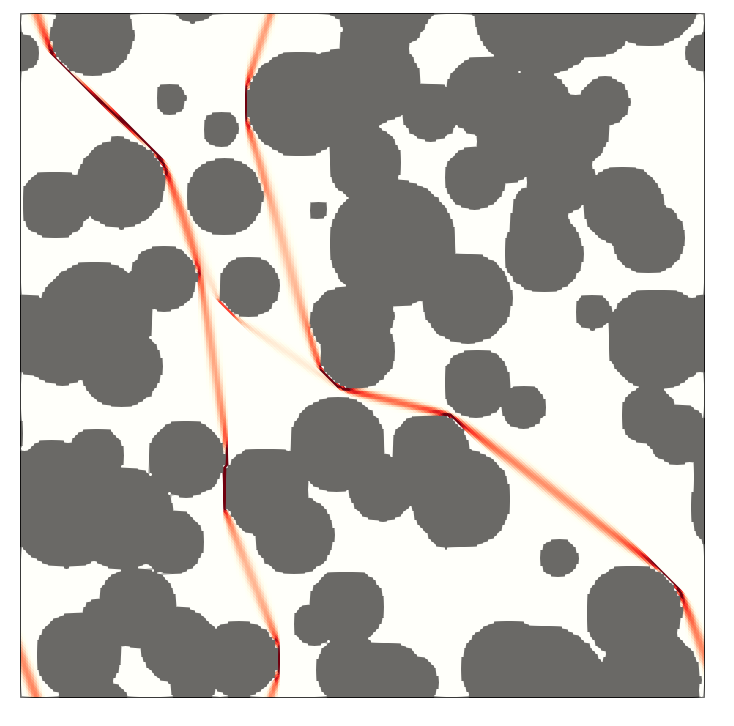}
\caption{$\theta=15^\circ$, $\Gceff=1.097\Gc$}
\end{subfigure}
\\
\begin{subfigure}{0.44\textwidth}
\includegraphics[width=\textwidth]{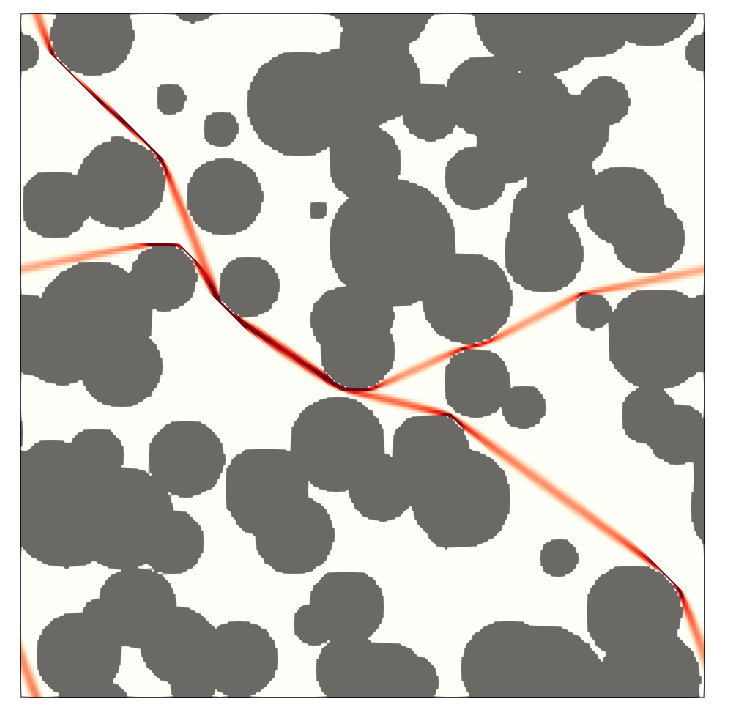}
\caption{$\theta=30^\circ$, $\Gceff=1.158\Gc$}
\end{subfigure}
\hfill
\begin{subfigure}{0.44\textwidth}
\includegraphics[width=\textwidth]{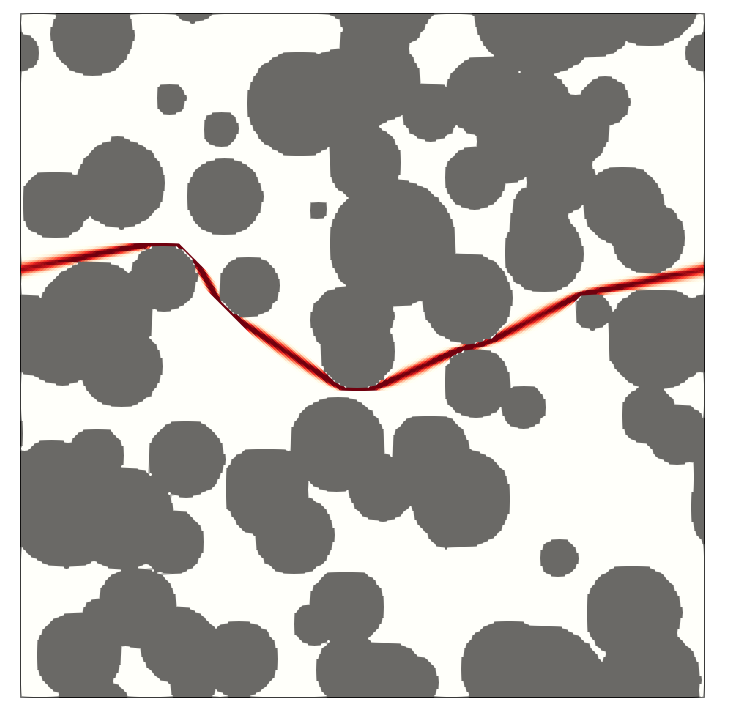}
\caption{$\theta=45^\circ$, $\Gceff=1.128\Gc$}
\end{subfigure}
\end{center}
\caption{Crack density field for a microstructure with random inclusions}\label{fig:minimal-crack-surface-random-inclusions-results}
\end{figure}

We now consider a similar setting in which the microstructure is obtained as a random collection of potentially overlapping disks of varying radius. We used a $250\times 250$ mesh. The obtained optimal crack surfaces have been represented in Figure \ref{fig:minimal-crack-surface-random-inclusions-results}. One can clearly see that purely vertical cracks are not allowed in this case for $\theta=0^\circ$, which results in an increased effective fracture energy by roughly 5\%. Similarly for $\theta=90^\circ$, the deviation from a horizontal crack is even more important, resulting in an increased fracture energy by 13\%. 

\section{Nonlinear membranes}\label{sec:membranes}
In this section, we consider the computation of the displacement field $\bu$ of an hyperelastic membrane subject to imposed tractions $\bT$ on the Neumann part of the boundary $\partial\Omega_\text{N}$ and to fixed displacement $\bu = 0$ on the Dirichlet part of the boundary $\partial \Omega_\text{D}$. In the finite-strain setting, the hyperelastic potential of a 3D material can be expressed as a function of some nonlinear strain measures. For instance, let us consider the Cauchy-Green strain tensor $\bC=\bF\T\bF$ where $\bF=\bI+\nabla \bu$ is the deformation gradient. The free energy hyperelastic potential is then $\psi(\bC)$ and the displacement field can be obtained as the solution to the following minimum principle:
\begin{equation}
\begin{array}{rl} \displaystyle{\inf_{\bu,\bC}} & \displaystyle{\int_\Omega \psi(\bC) \dOm - \int_{\partial \Omega_\text{N}} \bT\cdot\bu \dS}\\
\text{s.t.} & \bC = \bI + \nabla \bu + \nabla \bu\T +\nabla \bu\T\nabla \bu\\
& \bu = 0 \text{ on } \partial\Omega_\text{D}
\end{array}\label{eq:hyperelastic-min-principle}
\end{equation}

In the general case, such a problem is not convex which makes the analysis of hyperelastic materials rather complex. In particular, some equilibrium positions can be unstable and lead to buckling phenomena. 

\subsection{Tension field elastic membrane}
As regards thin hyperelastic membranes, local buckling (or wrinkling) will occur at very low load levels in compressed regions. In the limit of infinitely thin membranes, compressed stress states cannot be supported at all. Tension field theory \cite{wagner1929flat,reissner1938tension} has been proposed in order to simplify the analysis of thin membranes.

In the finite-deformation case, the tension-field theory has been first formalized by \cite{pipkin1994relaxed} by introducing a relaxed strain energy functional. More precisely, introducing the following quasi-convexification of $\psi$:
\begin{equation}
\begin{array}{rl}\psi_\text{memb}(\bC) =  \displaystyle{\inf_{\bCel}} & \psi(\bCel) \\
 \text{s.t. } &\bC_\text{el} \succeq \bC \end{array} \label{eq:relaxed-hyperelastic-potential}
 \end{equation}
the tension field variational principle is obtained when replacing $\psi$ with $\psi_\text{memb}$ in \eqref{eq:hyperelastic-min-principle}:
\begin{equation}
\begin{array}{rl} \displaystyle{\inf_{\bu,\bC}} & \displaystyle{\int_\Omega \psi_\text{memb}(\bCel) \dOm - \int_{\partial \Omega_\text{N}} \bT\cdot\bu \dS}\\
\text{s.t.} & \bCel \succeq \bI + \nabla \bu + \nabla \bu\T +\nabla \bu\T\nabla \bu\\
& \bu = 0 \text{ on } \partial\Omega_\text{D}
\end{array}\label{eq:relaxed-hyperelastic-min-principle}
\end{equation}

The above relaxed potential then provides a tension-field constitutive equation in terms of the second Piola-Kirchhoff stress $\bS$ as follows:
\begin{subequations}
\begin{align}
\bS &= 2\dfrac{\partial \psi}{\partial \bC}(\bC_\text{el})\\
\bC &= \bC_\text{el} + \bC_\text{w} \\
\bC_\text{w} &\preceq 0,\quad \bS \succeq 0, \quad \bS:\bC_\text{w}=0
\end{align}
\end{subequations}
where $\bC_\text{w}$ can be seen as an inelastic \textit{wrinkling} strain accounting for the occurrence of wrinkles in compressed regions. As a consequence, the resulting stress is always tensile.\\

\subsection{Conic reformulation}\label{sec:membrane-conic-reformulation}
Finally, Pipkin showed that when $\psi$ is a convex function of $\bC$, $\psi_\text{memb}$ turns out to be a convex function of $\bF$ or, equivalently, of the displacement gradient $\bG=\nabla \bu$ \cite{pipkin1994relaxed}. Indeed, if $\psi(\bCel)$ is convex, the relaxed minimum principle \eqref{eq:relaxed-hyperelastic-min-principle} is a convex program due to the following conic reformulation of the SDP constraint (see also \cite{kanno2011nonsmooth}):
\begin{equation}
\bCel \succeq \bC = \bI + \bG + \bG\T +\bG\T\bG
\end{equation}

Let us first recall the \textit{Schur complement lemma} for a PSD block-matrix:
\begin{lemma}
A symmetric block-matrix
\begin{equation}
\bZ = \begin{bmatrix}
\bU & \bV \\ \bV\T & \bW
\end{bmatrix}
\end{equation}
is positive semi-definite if and only if $\bW\succeq 0$ and $\bU-\bV\bW^{-1}\bV\T \succeq 0.$
\end{lemma}

Let us then consider the following symmetric matrix:
\begin{equation}
  \bZ = 
\begin{bmatrix}
\bCel & \bI+\bG\T \\
\bI+\bG & \bI
\end{bmatrix} \succeq 0
\end{equation}
we have $\bW=\bI \succ 0$ and $\bU-\bV\bW^{-1}\bV\T = \bCel - (\bI+\bG\T)(\bI+\bG) = \bCel - \bI - \bG - \bG\T - \bG\T\bG=\bC$. As a result, using the Schur complement lemma, $\bZ\succeq 0$ if and only if $\bCel \succeq \bC$. 

In conclusion, if $\psi(\bC)$ admits a convex conic representation in terms of $\bC$, problem \eqref{eq:relaxed-hyperelastic-min-principle} is a convex conic program.

\subsection{Ogden-type hyperelastic materials}
In \cite{kanno2011nonsmooth}, only Saint-Venant-Kirchhoff materials were considered for which $\psi(\bC)$ is a quadratic function. Unfortunately, this simple material model is not appropriate for modelling hyperelastic membranes. In the following, we consider an incompressible Ogden-type material (with only one term for simplicity) for which the 3D potential reads:
\begin{equation}
\psi(\lambda_1,\lambda_2,\lambda_3) = \dfrac{2\mu}{\alpha^2}\sum_{i=1}^3 \lambda_i^\alpha
\end{equation}
where $\mu$ is the shear modulus, $\alpha$ is some power-law exponent (which we assume in the following to be larger than 2) and $\lambda_1,\lambda_2,\lambda_3$ denote the principal stretches. The latter are the positive square roots of the Cauchy-Green tensor eigenvalues. Note that the case $\alpha=2$ corresponds to a neo-Hookean material.

In the incompressible case, we have $J=\lambda_1\lambda_2\lambda_3=1$ which yields a reduced energy:
\begin{equation}
\hat{\psi}(\lambda_1,\lambda_2)=\psi(\lambda_1,\lambda_2,\lambda_1^{-1}\lambda_2^{-1}) = \dfrac{2\mu}{\alpha^2}\left(\lambda_1^\alpha + \lambda_2^\alpha + \dfrac{1}{\lambda_1^\alpha\lambda_2^\alpha}\right)
\end{equation}

From the following relation for the principal stretches in 2D:
\begin{align}
\lambda_1^2 = \dfrac{1}{2}\left(C_{11}+C_{22} + \sqrt{(C_{11}-C_{22})^2 + 4C_{12}^2}\right)  \\
\lambda_2^2 = \dfrac{1}{2}\left(C_{11}+C_{22} - \sqrt{(C_{11}-C_{22})^2 + 4C_{12}^2}\right)  \\
\end{align}
one can see that the reduced energy can be equivalently expressed as a convex function of $\bC$ as follows:
\begin{equation}
\begin{array}{rl} \displaystyle{\hat{\psi}(\bC) = \min_{s,t}} & \dfrac{\mu}{2\beta^2}\left((t+s)^{\beta}+(t-s)^{\beta} + (t^2-s^2)^{-\beta}\right) \\
\text{s.t.} & t = \dfrac{1}{2}(C_{11}+C_{22}) \\
& \dfrac{1}{2}\sqrt{(C_{11}-C_{22})^2 + 4C_{12}^2} \leq s 
\end{array}
\end{equation}
where we assume that $\beta=\alpha/2 \geq 1$.

The above expression can further be re-expressed introducing additional auxiliary variables as:
\begin{equation}
\begin{array}{rl} \displaystyle{\hat{\psi}(\bC) = \min_{r,s,t,\by}} & \dfrac{\mu}{2\beta^2}(x_0+y_0+z_0) \\
\text{s.t.} & \begin{Bmatrix} (C_{11}-C_{22})/2 \\ C_{12} \end{Bmatrix} = \begin{Bmatrix} y_1 \\ y_2 \end{Bmatrix} \\
& \sqrt{y_1^2+y_2^2} \leq s \\
& (t+s)^\beta \leq x_0 \\
& (t-s)^\beta \leq y_0 \\
& (t^2-s^2) \geq u^2 \\
& u^{-2\beta} \leq z_0
\end{array}
\end{equation}
and where the last constraints can be formulated using a suitable quadratic and power cones as follows:
\begin{subequations}
\begin{align}
(t+s)^\beta \leq x_0 & \Leftrightarrow \begin{cases} x_2 = t+s \\ x_1 = 1 \\ |x_2|\leq x_0^{1/\beta} \text{ that is } \bx \in \Pp^{1/\beta}_3 \end{cases}\\
(t-s)^\beta \leq y_0 & \Leftrightarrow \begin{cases} y_2 = t-s \\ y_1 = 1 \\ |y_2|\leq y_0^{1/\beta} \text{ that is } \by \in \Pp^{1/\beta}_3 \end{cases}\\
(t^2-s^2) \geq u^2 & \Leftrightarrow \sqrt{s^2+u^2}\leq t \text{ that is } (t,s,u)\in \Qq_3\\
u^{-2\beta}\leq z_0 & \Leftrightarrow 1 \leq  u z_0^{2\beta} \Leftrightarrow 1\leq u^{\frac{1}{1+2\beta}} z_0^{\frac{2\beta}{1+2\beta}} \text{ that is } (u,z_0,1) \in \Pp_3^{\frac{1}{1+2\beta}}
\end{align}
\end{subequations}

\subsection{Material point response validation}
We formulate the variational problem \eqref{eq:relaxed-hyperelastic-min-principle} using the incompressible Ogden-type behaviour expressed using the previous conic constraints. We simulate the response of a material point by prescribing an affine displacement field on the boundary:
\begin{equation}
\bu = \overline{\bG}\cdot\bx \quad \text{on }\partial\Omega 
\end{equation}
with an imposed displacement gradient of the form:
\begin{equation}
\overline{\bG} = \begin{bmatrix}
\delta_x & \gamma \\ 0 & \delta_y
\end{bmatrix} \Rightarrow \bC = \begin{bmatrix}
(1+\delta_x)^2 & \gamma(1+\delta_x) \\ \gamma(1+\delta_x) & (1+\delta_y)^2+\gamma^2
\end{bmatrix}
\end{equation}
where $\delta_x$ (resp. $\delta_y$) is a horizontal (resp. vertical) elongation and $\gamma$ a shear distortion. These specific boundary conditions induce a uniform deformation state that is $\bG=\nabla \bu = \overline{\bG}$ inside $\Omega$.
The principal true (or Cauchy) stresses in the incompressible case for a plane-stress membrane ($\sigma_3=0$) are then given by:
\begin{equation}
\sigma_i = \lambda_i\dfrac{\partial \hat\psi}{\partial \lambda_i} = \dfrac{2\mu}{\alpha}\left(\lambda_i^\alpha - \dfrac{1}{\lambda_1^\alpha\lambda_2^\alpha}\right), \quad i=1,2 \label{eq:membrane-stress-theoretical}
\end{equation}

In our numerical implementation, the stress is recovered from the Lagrange multiplier associated with the SDP constraint $\bZ\succeq 0$. Indeed, one can easily show from SDP duality that the latter dual variable $\bLambda_{\bZ}$  is a block-matrix which is also positive:
\begin{equation}
\bLambda_{\bZ} = \begin{bmatrix}
\bLambda_{\bU} & \bLambda_{\bV} \\ \bLambda_{\bV}\T & \bLambda_{\bW}
\end{bmatrix} \succeq 0
\end{equation}
and that one has at the optimum:
\begin{equation}
\bLambda_{\bU} = \dfrac{\partial \psi}{\partial \bC}(\bCel) \succeq 0
\end{equation}
which gives the Piola-Kirchhoff stress $\bS$ up to a factor 2.

\begin{figure}
\begin{center}
\begin{subfigure}{0.49\textwidth}
\includegraphics[width=\textwidth]{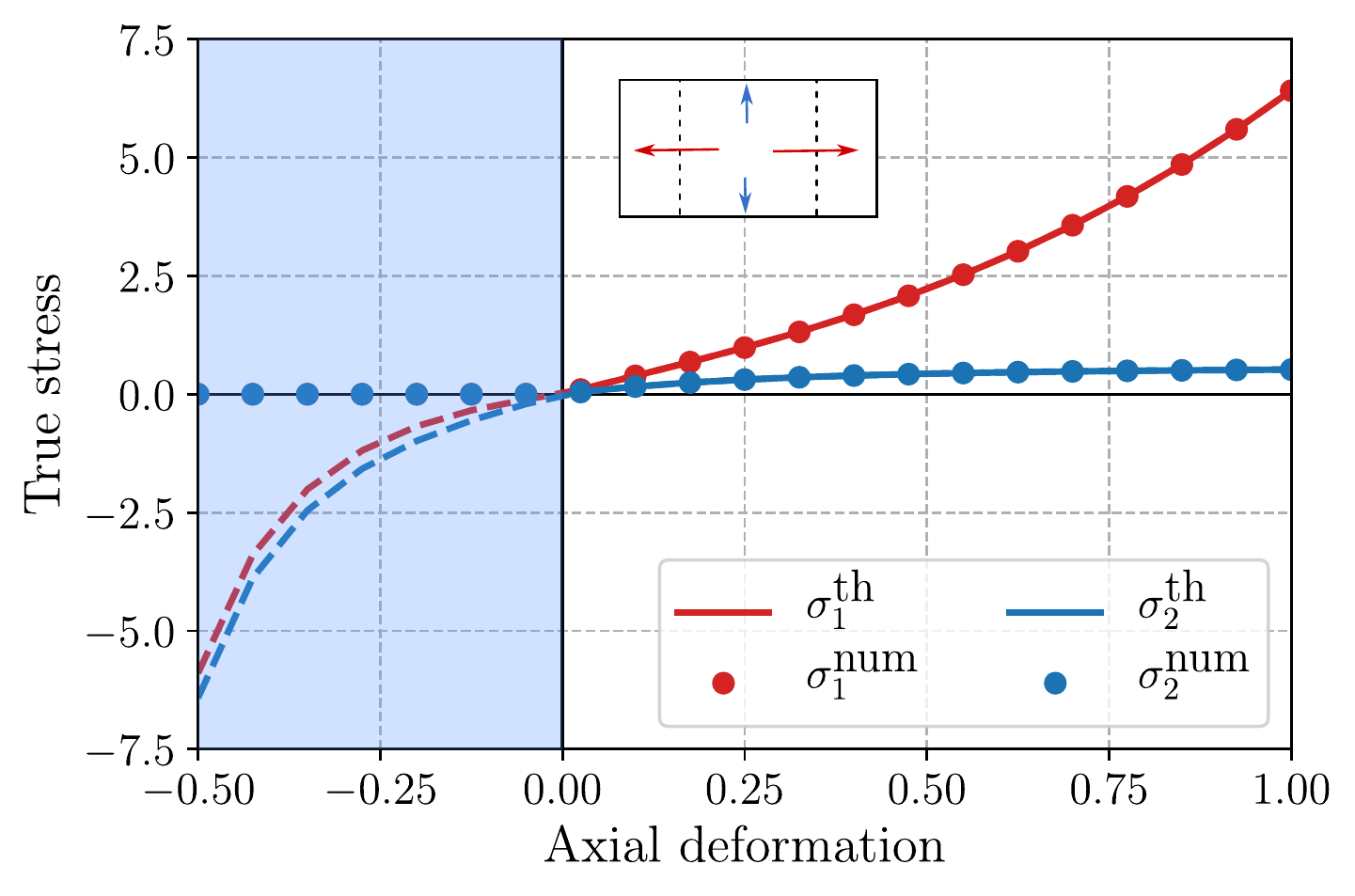}
\caption{$\varepsilon_y=\delta=0$}\label{fig:membrane-stress-uniaxial}
\end{subfigure}
\hfill
\begin{subfigure}{0.49\textwidth}
\includegraphics[width=\textwidth]{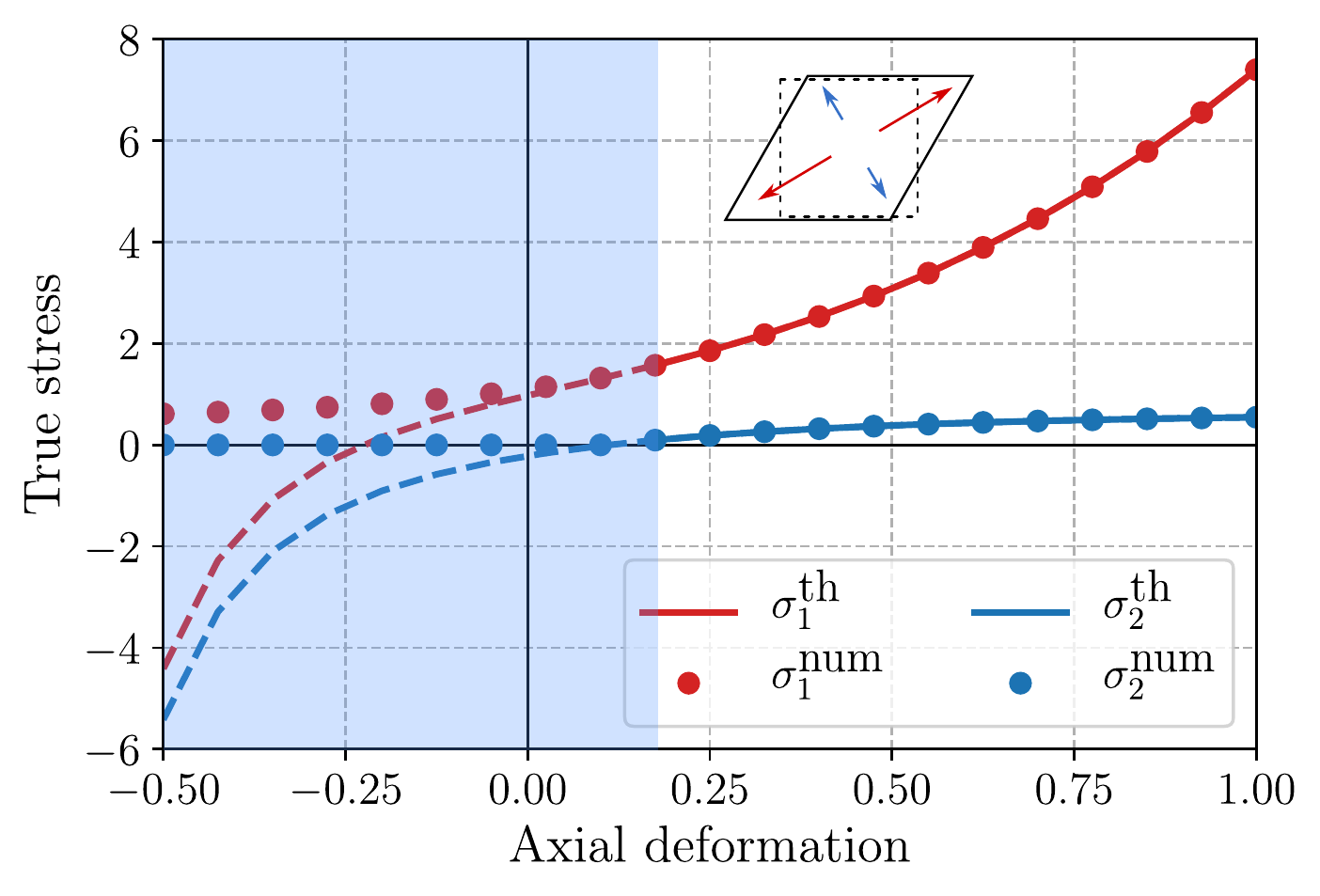}
\caption{$\varepsilon_y=0.05$, $\delta=0.5$}\label{fig:membrane-stress-combined}
\end{subfigure}
\end{center}
\caption{Symbols correspond to the numerical tension field state and solid/dashed lines correspond to expression \eqref{eq:membrane-stress-theoretical}. The latter is not valid in the light blue region (dashed lines) which denotes the wrinkled state where $\sigma_2 = 0$.}\label{fig:membrane-material-point}
\end{figure}

In Figure \ref{fig:membrane-material-point}, we report results for the evolution of both principal true stresses as a function of an imposed axial elongation $\varepsilon_x$ either in a pure uniaxial elongation case $\varepsilon_y=\delta=0$ (Fig. \ref{fig:membrane-stress-uniaxial}) or in a combined elongation and shear case $\varepsilon_y=0.05$, $\delta=0.5$ (Fig. \ref{fig:membrane-stress-combined}). One can see that the numerical tension field results match exactly with the theoretical expression \eqref{eq:membrane-stress-theoretical} in regions where both $\sigma_1,\sigma_2 > 0$ (solid lines). In the remaining region where the analytical $\sigma_2 \leq 0$, the numerical tension field stresses deviate from equation \eqref{eq:membrane-stress-theoretical} (dashed lines) and yields a pure tension state with $\sigma_2=0$ and $\sigma_1>0$ which corresponds to wrinkles having formed along principal direction 2.

\subsection{Annular square membrane deformation}

We illustrate our implementation on a square annular membrane initially located in the $(Ox_1x_2)$ plane, embedded in $\RR^3$. Its outer boundary of size $W_\text{out}=50$ mm is fixed whereas the inner boundary, of size $W_\text{in}=17.5$ mm, is subjected to an in-plane torsion of angle $90^\circ$ and to an out-of-plane vertical displacement of amplitude $t W_\text{out}/2$ for $t=0$ to $t=1$. Note that a similar example has been investigated in \cite{de2015finite} with a Neo-Hookean model whereas we use here $\alpha=3.5$.

In this setting, since we consider out-of-plane displacements, the displacement and deformation gradients $\bG$ and $\bF$ are in fact $3\times 2$ matrices such that:
\begin{equation}
\bF = \overline{\bI} + \bG\quad \text{where } \overline{\bI} = \begin{bmatrix}
\be_1 & \be_2
\end{bmatrix}\in \mathbb{R}^{3\times 2}
\end{equation}
where $\be_1$, $\be_2$ denotes orthonormal basis vectors of the membrane initial reference plane. The strain tensors $\bC,\bCel$ are still $2\times 2$ matrices. The only modification to apply to section \ref{sec:membrane-conic-reformulation} concerns the SDP constraint. We now have:
\begin{equation}
\bCel \succeq \bI_2 + \overline{\bI}\T\bG + \bG\T\overline{\bI} + \bG\T\bG
\end{equation}
where $\bI_n$ is the identity matrix of dimension $n$. The latter constraint can now be equivalently reformulated as follows:
\begin{equation}
\bZ = \begin{bmatrix}
\bCel & \overline{\bI}\T+\bG\T \\
\overline{\bI} + \bG & \bI_3
\end{bmatrix} \succeq 0
\end{equation}

\begin{figure}
\begin{center}
\begin{subfigure}{0.3\textwidth}
\includegraphics[width=\textwidth]{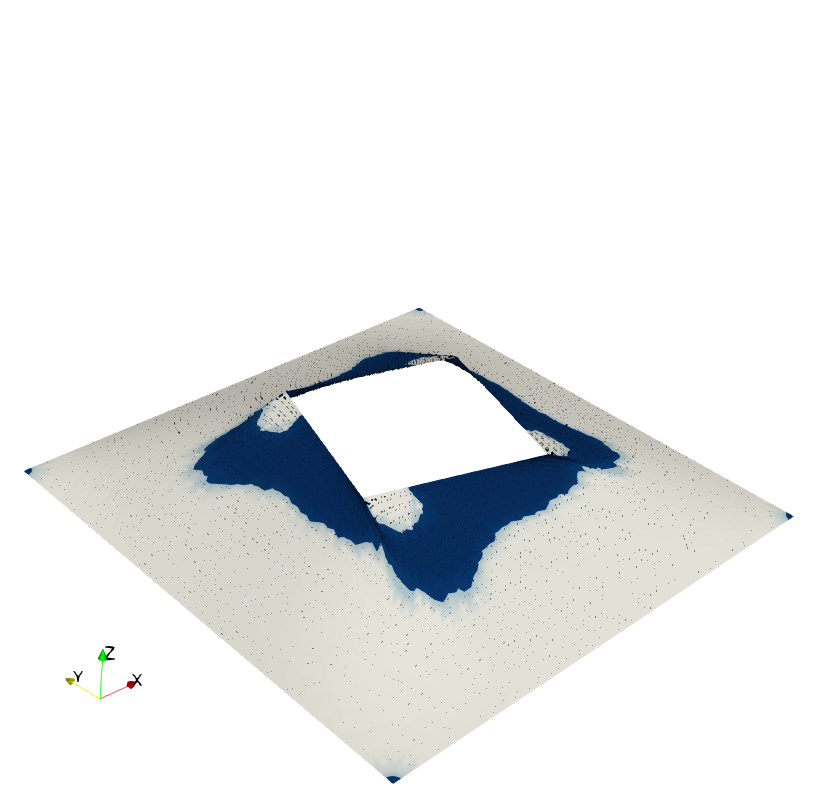}
\caption{$t=0.2$}
\end{subfigure}
\hfill
\begin{subfigure}{0.3\textwidth}
\includegraphics[width=\textwidth]{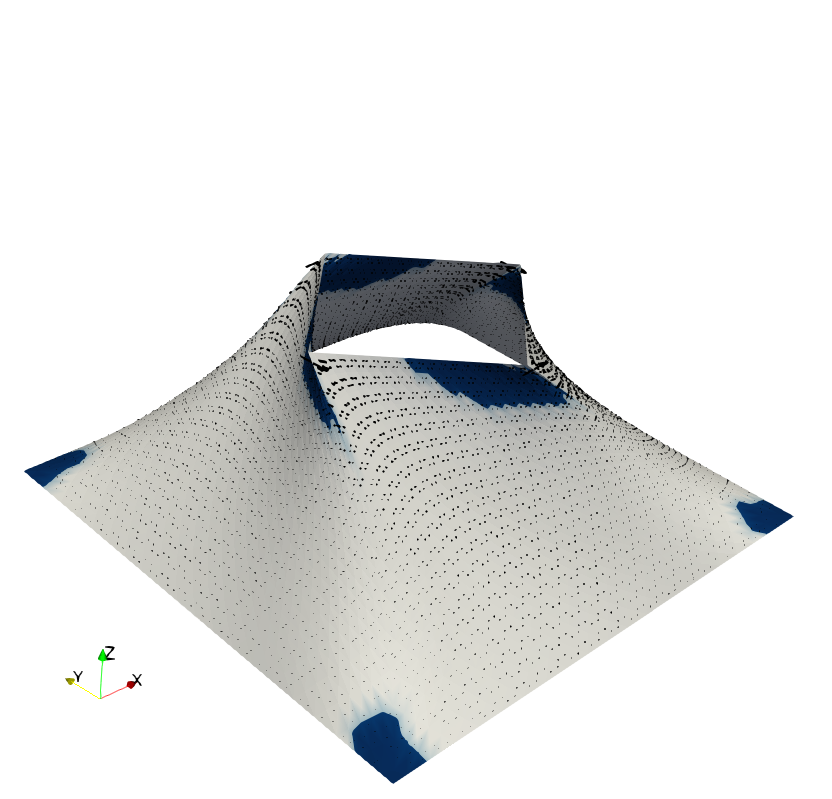}
\caption{$t=0.5$}
\end{subfigure}
\hfill
\begin{subfigure}{0.3\textwidth}
\includegraphics[width=\textwidth]{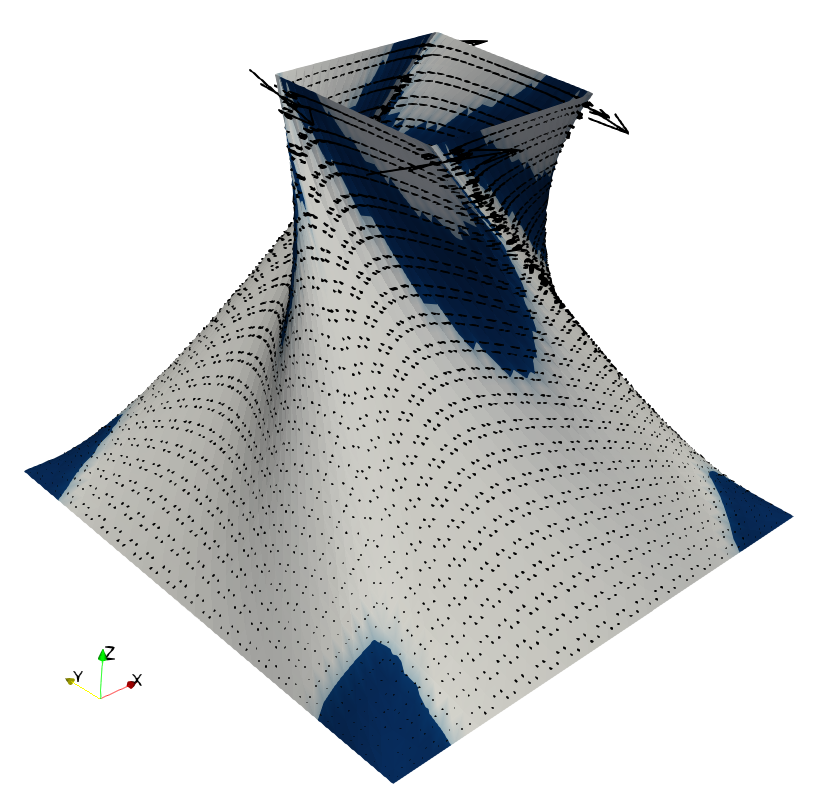}
\caption{$t=1$}
\end{subfigure}
\end{center}
\caption{Annular square deformed membrane at various loading steps. Blue regions denote wrinkled regions. Black arrows indicate principal stresses amplitude and orientation.}\label{fig:membrane-annular}
\end{figure}

Figure \ref{fig:membrane-annular} reports different snapshots of the deformation and wrinkled state of the membrane. We can observe that initially a large inner region is in a wrinkled state due to compression induced by the torsion. At larger load steps, the wrinkled region extent tends to diminish due to the important tension exerted by the vertical displacement. In the last stage, the shape of the wrinkled regions is quite complex due to the combination of excessive shearing and elongation. The obtained wrinkled regions are quite similar to those obtained in \cite{de2015finite} with a slightly different material model. Finally, let us remark that it was not necessary to subdivide the final loading into smaller load steps. Interior-point methods are known to be remarkably robust and quite insensitive to a good or bad initial guesses for the solution. The final solution, as well as even more extremely deformed configurations for $t>1$, could be obtained in only 20 iterations.

\section{Viscoplastic fluid flows in Hele-Shaw cells}\label{sec:viscoplastic}
\begin{figure}
\begin{center}
\includegraphics[width=0.7\textwidth]{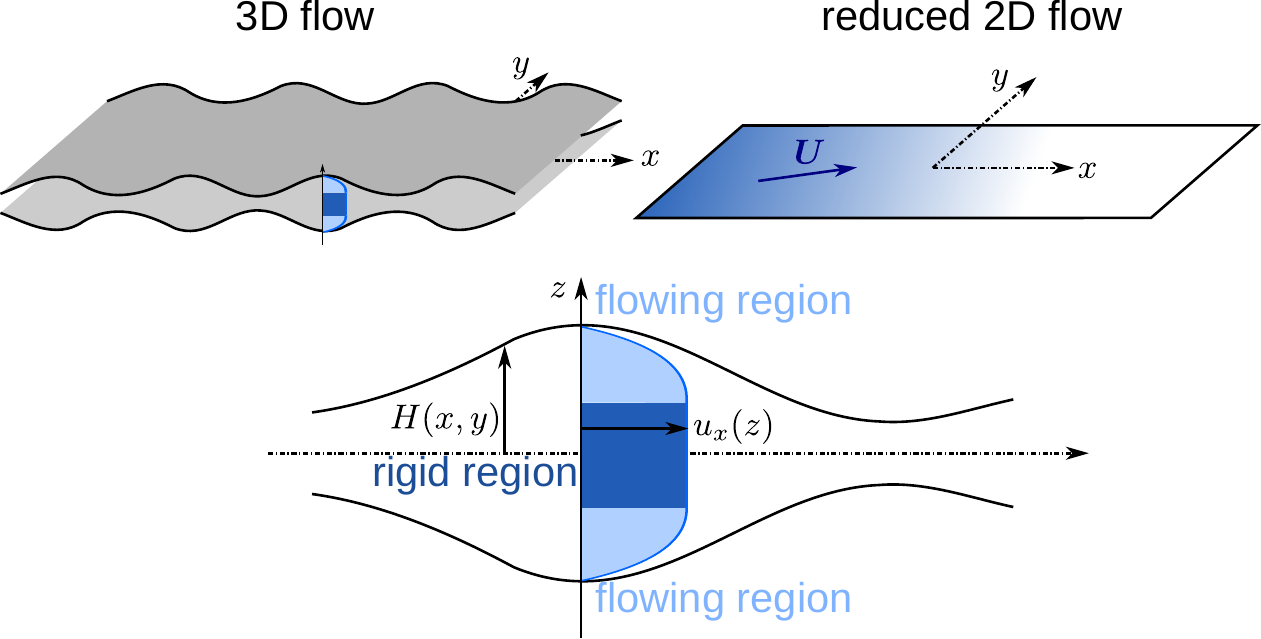}
\end{center}
\caption{Flow in a Hele-Shaw cell for a viscoplastic fluid and dimensional reduction}
\end{figure}

A Hele-Shaw cell consists of two plates which are extremely close to each other in order to induce a bi-dimensional flow. This serves as a model for flows in porous media. Each plate is located at a position $z=\pm H(x,y)$ (symmetric configuration with respect to $z=0$).

Although the flow is macroscopically two-dimensional and can be described by a macroscopic in-plane velocity $\bU = U_x(x,y)\be_x+U_y(x,y)$, the local three-dimensional flow $\bu(x,y,z)$ admits a specific profile in the transverse $z$-direction between both plates. 

In order to link the local 3D behaviour of the fluid with an effective 2D behaviour, a dimensional reduction procedure is considered accounting for several hypotheses \cite{huilgol2015fluid}:
\begin{itemize}
\item the transverse velocity is negligible: $u_z\approx 0$
\item the velocity variations in the transverse directions are much larger than the in-plane variations $\|\bu_{,x}\|,\|\bu_{,y}\| \ll \|\bu_{,z}\|$
\item a no-slip condition is assumed along both plates walls $z=\pm H$
\item inertia and body forces can be ignored
\end{itemize}

Such a procedure is not trivial in the case of a Bingham fluid \cite{dean2007numerical,bleyer2018advances} which exhibits solid and flowing regions along the transverse direction. The yield point between flowing and solid region is unknown a priori and depends upon the corresponding stress and pressure gradient.

Indeed, with the above hypotheses, one can easily show that the shear stress field varies linearly along the transverse direction as follows:
\begin{align}
\bsig(x,y,z) &\approx \begin{bmatrix} 0 & 0 & \tau_x \\ 0 & 0 & \tau_y \\ \tau_x & \tau_y & 0\end{bmatrix} \\
\text{with } \btau(x,y) &= z \nabla p(x,y)
\end{align}
where $\btau(x,y,z)$ is the anti-plane shear stress vector and $p(x,y)$ is the fluid pressure.

\subsection{Determining effective potentials}

The Hele-Shaw effective behaviour can be described through effective potentials, either in stress-based form as a potential function $\Psi(\bG)$ of the pressure gradient $\bG(x,y)=-\nabla p(x,y)$ or in velocity-based form as a potential function $\Phi(\bU)$ of the macroscopic velocity $\bU(x,y)$. Both potentials are dual to each other, that is they are related via Legendre-Fenchel transform $\Psi = \Phi^*$ and we have:
\begin{equation}
\bU = \partial_{\bG} \Psi(\bG) \quad,\quad \bG = \partial_{\bU} \Phi(\bU)
\end{equation}

The stress-based potential for instance can be obtained from the integrated stress potential as follows:
\begin{equation}
\Psi(\bG) = \dfrac{1}{2H}\int_{-H}^H \psi(\bsig(x,y,z))\dz
\end{equation}
where $\psi$ is the 3D stress-based potential of the constituting fluid. 

For instance, for a Newtonian fluid, one has $\psi(\bsig) = (\bsig:\bsig)/(4\mu)$ with $\bd = \partial_{\bsig} \psi = \bsig/(2\mu)$. The corresponding effective potential is then given by:
\begin{equation}
\Psi(\bG) = \int_{-H}^H \dfrac{1}{4\mu H} z^2 \bG\cdot\bG \dz = \dfrac{H^2}{6\mu} \bG\cdot\bG
\end{equation}
which results in the following linear effective constitutive relation:
\begin{equation}
\bU = \partial_{\bG} \Psi(\bG) = \dfrac{H^2}{3\mu} \bG
\end{equation}
from which we recover the classical Darcy equation between two parallel plates.

Alternatively, the velocity-based potential is obtained via Legendre-Fenchel transform:
\begin{equation}
\Phi(\bU) = \sup_{\bG} \{\bU\cdot\bG - \Psi(\bG)\}
\end{equation}
which after standard convex duality computations yields:
\begin{equation}
\begin{array}{rl} \displaystyle{\Phi(\bU) = \inf_{\bgamma(z)}} & \displaystyle{\dfrac{1}{2H}\int_{-H}^H \phi(\bd(z))\dz} \\
\text{s.t.} & \bd(z) = \begin{bmatrix} 0 & 0 & \gamma_x(z) \\ 0 & 0 & \gamma_y(z) \\ \gamma_x(z) & \gamma_y(z) & 0 \end{bmatrix} \\
& \bU + \displaystyle{\dfrac{1}{2H}\int_{-H}^{H} z\bgamma(z)\dz = 0}
\end{array}
\end{equation}

In the above formulation, we see that when interpreting $\bgamma(z)$ as the local strain rate $\bgamma(z)=\bu_{,z}$, the constraint indeed enforces the link with the macroscopic velocity:
\begin{equation}
\bU + \dfrac{1}{2H}\int_{-H}^{H} z\bu_{,z}\dz = \bU - \dfrac{1}{2H}\int_{-H}^{H} \bu\dz =0
\end{equation}

\subsection{The Bingham case and its conic representation}

\subsubsection{Numerical approximation}
For a Bingham fluid of viscosity $\mu$ and yield stress $\tau_0$, the velocity-based potential can be expressed as:
\begin{equation}
\begin{array}{rl} \displaystyle{\Phi(\bU) = \inf_{\bgamma(z)}} & \displaystyle{\dfrac{1}{2H}\int_{-H}^H \left(\dfrac{\mu}{2}\|\bgamma(z)\|^2 + \tau_0\|\bgamma(z)\|\right)\dz} \\
\text{s.t.} &  \displaystyle{\bU + \dfrac{1}{2H}\int_{-H}^{H} z\bgamma(z)\dz = 0}
\end{array}
\end{equation}

Unfortunately, there is no closed-form expression for solving the inner minimization over the $\bgamma(z)$ field. For a practical numerical implementation, we therefore have to resort to some numerical approximation. In the following, we choose to restrict the field $\bgamma(z)$ to a finite set of values $\bgamma_i = \bgamma(z_i)$ at some integration points $z_i$. The above integrals are then replaced by a numerical quadrature using such points. Moreover, we benefit from our initial assumption of a symmetric configuration so that the integral is performed over $[0;H]$ only. Introducing, the non-dimensional integration points $\xi_i=z_i/H$ and their corresponding quadrature weights $\omega_i$ on $[0;1]$, we finally have the following approximate potential depending on the number $m$ of quadrature points:

\begin{equation}
\begin{array}{rl} \displaystyle{\Phi_m(\bU) = \inf_{\bgamma_i}} & \displaystyle{\sum_{i=1}^m \omega_i\left(\dfrac{\mu}{2}\|\bgamma_i\|^2 + \tau_0\|\bgamma_i\|\right)} \\
\text{s.t.} & \displaystyle{\bU + H \sum_{i=1}^m \omega_i\xi_i\bgamma_i = 0}
\end{array}
\end{equation}

Note that such an approach has already been proposed in \cite{bleyer2016numerical} when deriving effective strength criteria of shells from the local 3D strength criterion.

\subsubsection{Definition using the generalized marginal operator}

Implementing directly the above representation lacks generality since it must be repeated each time one wants to change the underlying fluid potential e.g. considering a Hershel-Bulkley fluid for instance. Instead, we can use the generalized "marginal" operator \eqref{eq:marginals} so that any Hele-Shaw potential $\Psi$ can be obtained as the "marginal" of the underlying fluid potential $\phi(\bgamma)$ through the linear operators $\bA_i = -H\omega_i \xi_i \bI$ and the coefficients $c_i=\omega_i$. Note that the Lagrange multiplier associated with the last constraint gives directly access to the pressure gradient $\bG=\partial_{\bU} \Phi(\bU)$.

\subsection{Validation of the quadrature discretization}
\begin{figure}
\begin{center}
\includegraphics[width=0.8\textwidth]{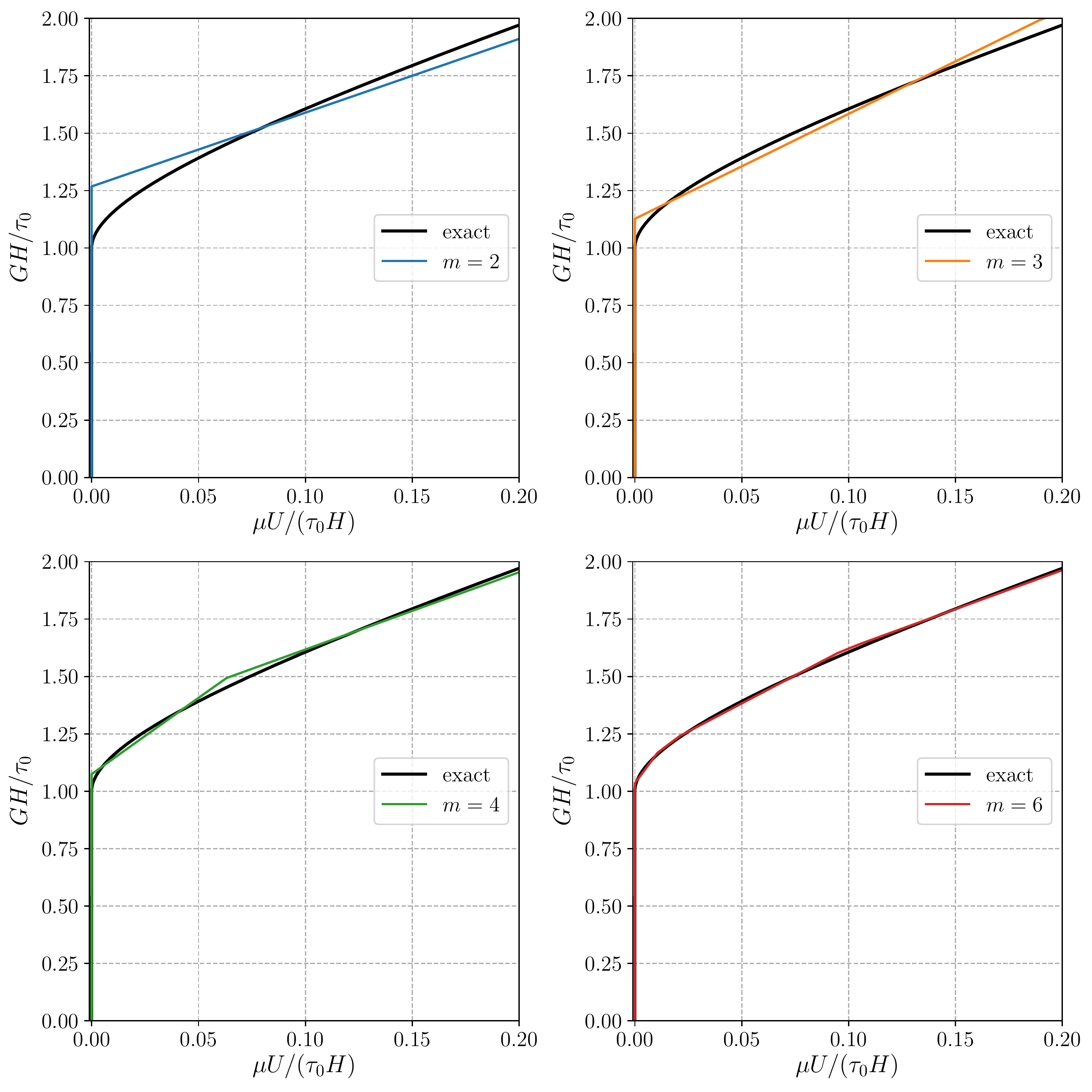}
\end{center}
\caption{Flow curves (norm of the pressure gradient $G=\|\bG\|$ as a function of filtration velocity magnitude $U=\|\bU\|$) of the Bingham Hele-Shaw model for different quadrature points $m$. Black curves correspond to \eqref{eq:hele-shaw-exact}.}\label{fig:hele-shaw-material-point}
\end{figure}

We now compare in Figure \ref{fig:hele-shaw-material-point} the obtained flow curves with the analytical solution provided in \cite{huilgol2015fluid}:
\begin{align}
\bU &= \dfrac{H^2}{3\mu}\left(1 - \dfrac{3}{2\|\bG\|} + \dfrac{1}{2\|\bG\|^3} \right)\bG \quad \text{ if }\|\bG\|>\tau_0/H \label{eq:hele-shaw-exact}\\
\bU &= 0 \text{ otherwise} \notag 
\end{align}

One can see that with $m=2$ quadrature points, the behaviour resembles that of a Bingham fluid (piecewise linear flow curve) with an overestimated yield stress and an underestimated asymptotic regime. Introducing more points enables to approach the non-linear black curve with a piecewise linear curve of better quality. The choice $m=4$ in particular is already extremely satisfying.

\subsection{Flow in a random porous medium}

\subsubsection{Flow equations and variational principle}

We now consider the Hele-Shaw flow on some domain $\Omega$. The main unknowns for this problem are the filtration velocity $\bU$ and the pressure field $p$. They must satisfy the balance equation and mass conservation condition:
\begin{align}
\nabla p &= 0 \quad \text{in } \Omega \\
\div\bU &= 0 \quad \text{in } \Omega
\end{align}
and they are related through the Hele-Shaw potential:
\begin{equation}
\bG = -\nabla p = \partial_{\bU} \Phi(\bU)
\end{equation}
along with some boundary conditions on pressure or velocity:
\begin{align}
p &= p_0 &\text{ on } \partial\Omega_\text{D} \\
\bU\cdot \bn &= q &\text{ on } \partial\Omega_\text{N}
\end{align}
where $\bn$ is the outward unit normal, $p_0$ and imposed pressure and $q$ an imposed flux.\\

One can then easily see that the above equations can be derived from the following filtration velocity variational principle:
\begin{equation}
\begin{array}{rl}
\displaystyle{\min_{\bU}} &\displaystyle{ \int_{\Omega}\Phi(\bU)\text{d}\Omega - \int_{\partial\Omega_\text{D}} p_0 \bU\cdot\bn \dS}\\
\text{s.t.} & \div\bU = 0  \text{ in }\Omega\\
& \bU\cdot\bn = q\text{ on } \partial\Omega_\text{N}
\end{array} \label{eq:hele-shaw-var-principle}
\end{equation}
with the pressure field being obtained from the Lagrange multiplier associated with the mass conservation equation.

\subsubsection{A random medium setting}

As an illustration, we aim at simulating the flow of such a fluid in a random porous-like medium. In particular, we consider a medium with a varying heights distribution $H(x,y)$ for which some regions are almost completely obstructed ($H(x,y)\approx 0$). In such situations, it is expected to find channelized flows near the yield point when macroscopic flow starts to occur, as discussed in details in \cite{hewitt2016obstructed} for instance.\\

We consider here a unit square domain $\Omega = [0;1]\times[0;1]$, subjected to $\bU\cdot \bn=0$ on the bottom and top surfaces and imposed pressure $p=0$ on the left boundary and $p=p_\text{out}$ on the right boundary, mimicking an imposed macroscopic pressure gradient $\overline{G} = p_\text{out}/L$ with $L=1$ here. In the following, we choose $\mu=1$ and $\tau_0=1$ as well as $m=4$ for the Hele-Shaw potential representation.

\begin{figure}
\begin{center}
\includegraphics[width=0.5\textwidth]{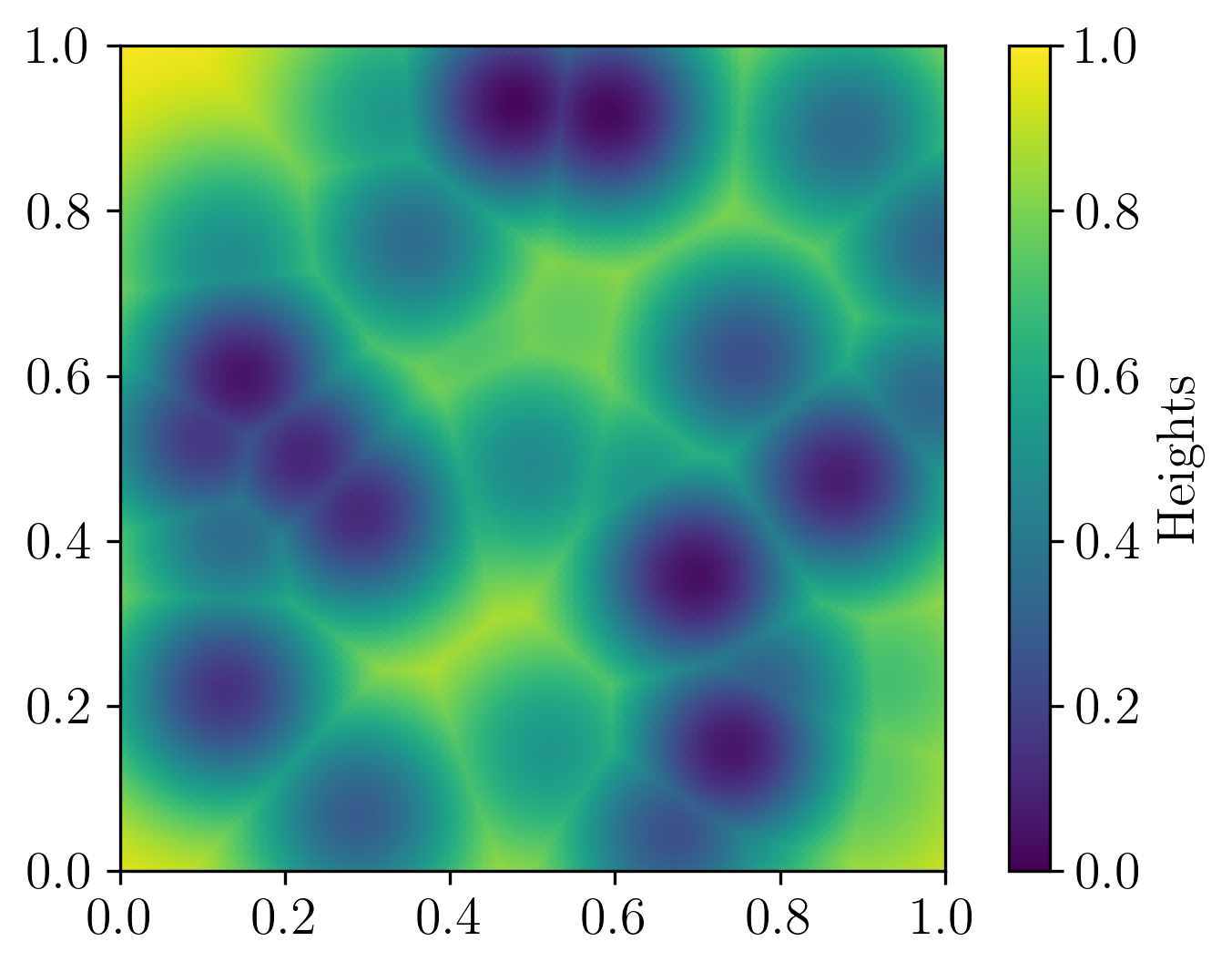}
\end{center}
\caption{Random realization of a Hele-Shaw cell with spatially varying height.}\label{fig:hele-shaw-heights}
\end{figure}

The Hele-Shaw cell consists of a heights distribution $H(x,y)\in[0;1]$ which is randomly generated through the piecewise combination of Gaussian profiles centered around random seed points (see Figure \ref{fig:hele-shaw-heights}). The field of heights is represented on a piecewise constant space which will also be used to discretize the pressure field. The velocity field is discretized using a Brezzi-Douglas-Marini space of degree $k=1$. The pressure field being a discontinuous function of degree $k-1=0$, this pair of function spaces is known to provide a stable approximation for such mixed Poisson-like systems \cite{brezzi1985two}. The mass conservation equation $\operatorname{div}(\bU)=0$ is enforced weakly using the pressure as a Lagrange multiplier.\\

\begin{figure}
\begin{center}
\begin{subfigure}{0.32\textwidth}
\includegraphics[width=\textwidth]{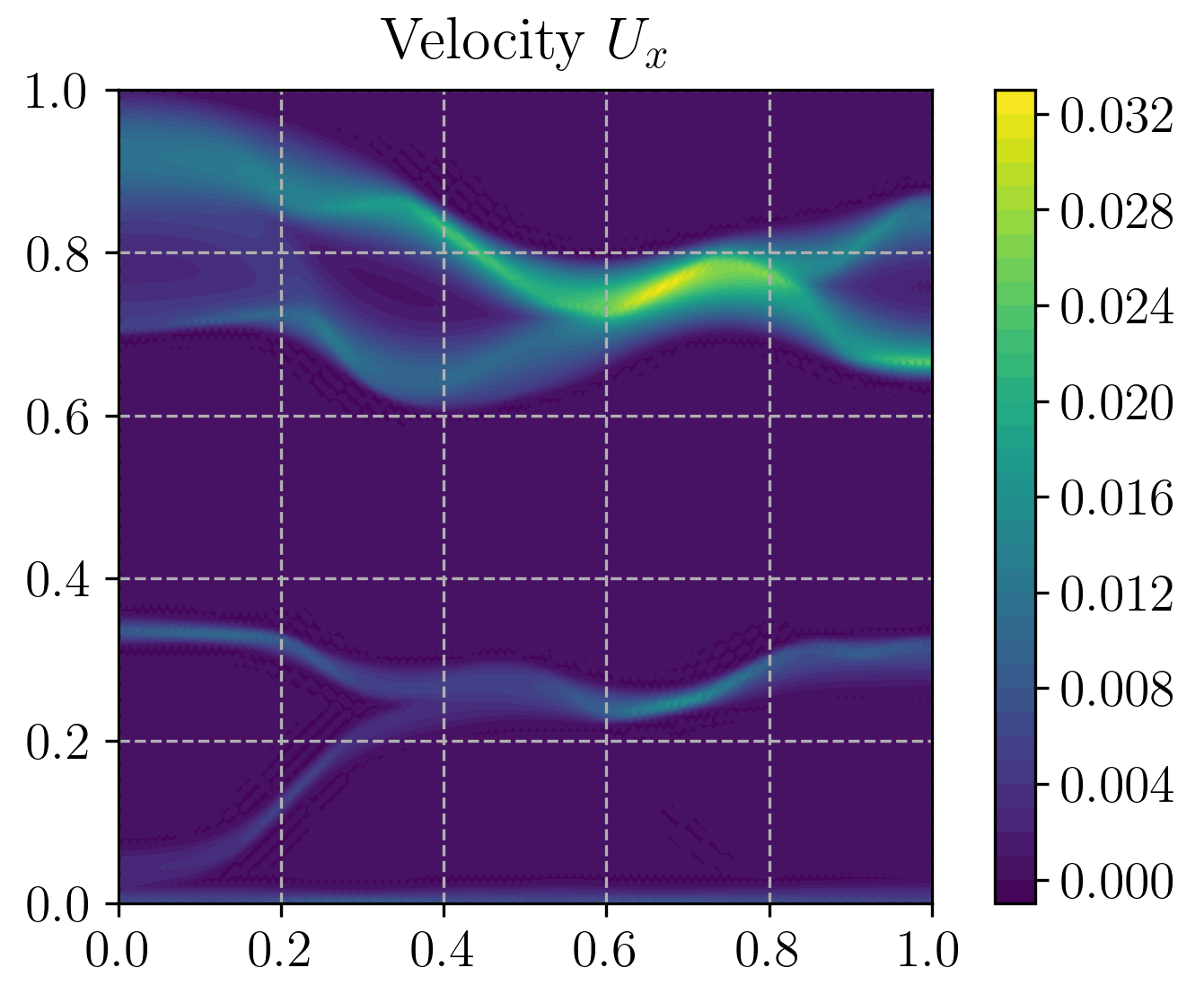}
\caption{$\overline{G}=2$}\label{fig:hele-shaw-flow-p-2}
\end{subfigure}
\hfill
\begin{subfigure}{0.32\textwidth}
\includegraphics[width=\textwidth]{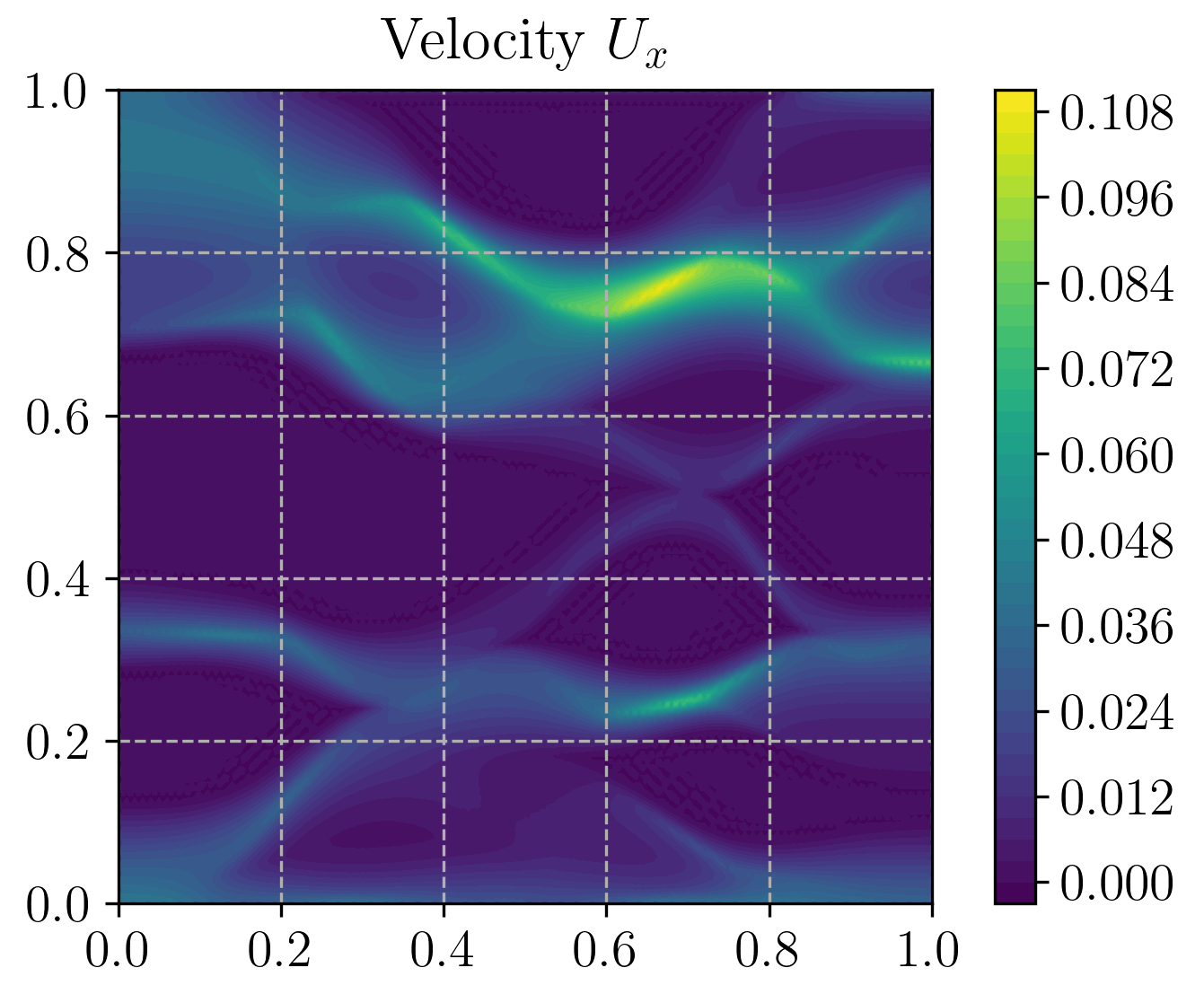}
\caption{$\overline{G}=2.5$}\label{fig:hele-shaw-flow-p-2.5}
\end{subfigure}
\hfill
\begin{subfigure}{0.32\textwidth}
\includegraphics[width=\textwidth]{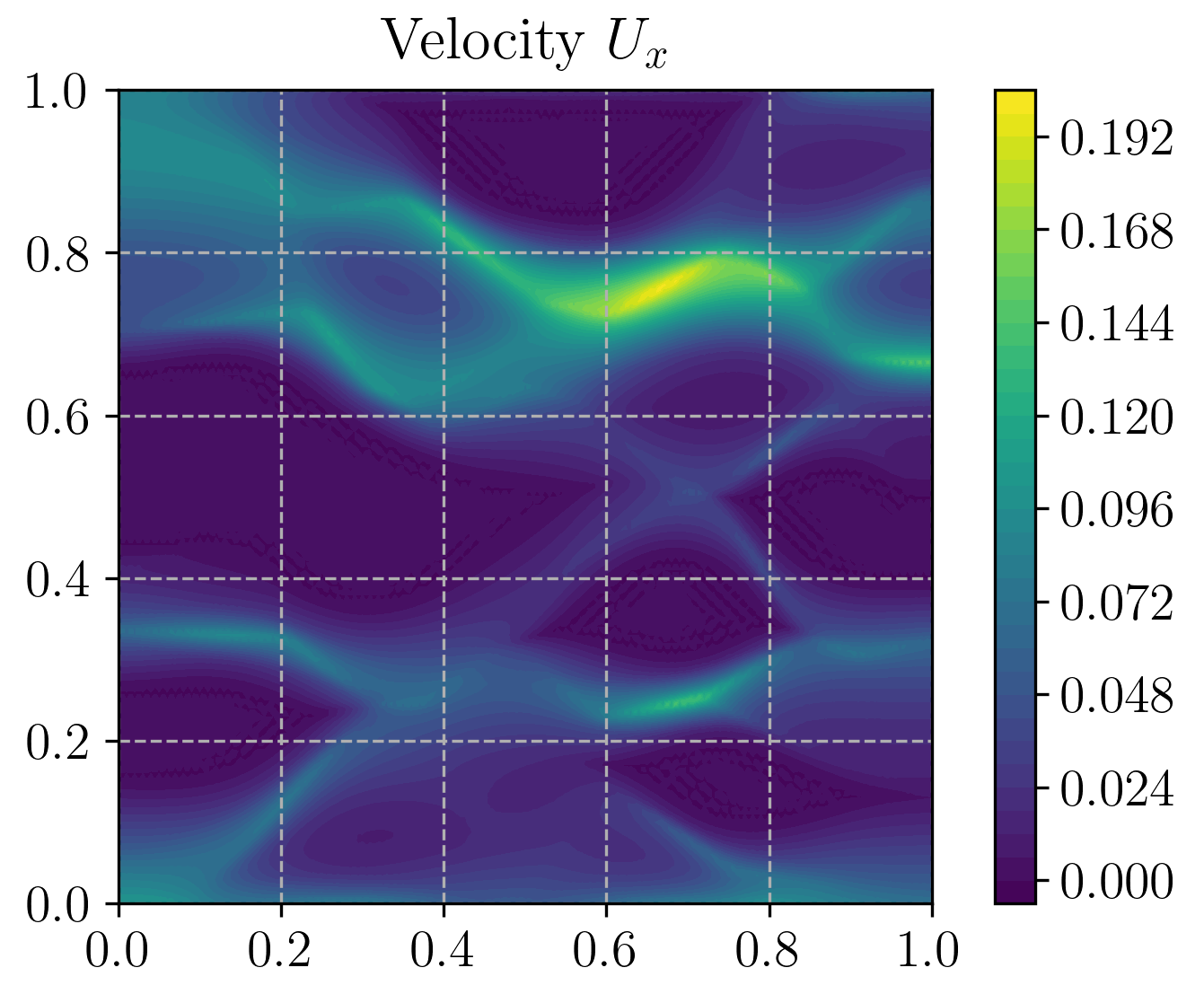}
\caption{$\overline{G}=3$}\label{fig:hele-shaw-flow-p-3}
\end{subfigure}
\end{center}
\caption{Horizontal filtration velocity maps for different imposed pressure gradients $\overline{G}$.}
\end{figure}

Solving the Hele-Shaw variational principle \eqref{eq:hele-shaw-var-principle} indeed results in the following observations. First, the local velocity is always zero if the applied pressure is below a certain threshold (here around $\overline{G} \approx 1.5$) and above this threshold, the flow starts occurring in localized channels only (see Figure \ref{fig:hele-shaw-flow-p-2}). For higher pressure gradients (Figure \ref{fig:hele-shaw-flow-p-3}), the flow then resembles that of a classical Newtonian fluid. 

\section{Conclusions and outlook}\label{sec:conclusions}

In this work, we showed how conic programming can be applied to formulate various complex problems arising in non-smooth mechanics. We have discussed how the incremental variational problem of standard generalized materials can be formulated as a conic program provided both the free energy density and the dissipation pseudo-potential are convex functions which admit a conic representation. However, such a satisfying convex variational framework is not always applicable. For instance, softening behaviours introduce formidable challenges as the corresponding potentials become non-convex, resulting in mathematical ill-posedness, existence of localized solutions and mesh-dependency of the corresponding numerical methods. Similarly, non-associated behaviours depart from the existence of a potential but can be formulated as mixed complementary programs. Alternative strategies relying on iterative resolutions of classical conic programs have already been proposed \cite{portioli2014limit,nodargi2021limit} and should deserve more investigation in the future.

Finally, convex programs also frequently appear in the context of optimal design and topology optimization. In this context, uncertainties on problem data are often difficult to include and may strongly influence the optimal solution. Robust optimization \cite{ben2009robust}, which often benefits from a tractable conic programming reformulation, might be an interesting approach for taking into account uncertainties in a conic programming framework.

\section*{Acknowledgements}
The author acknowledges Alexander Niewiarowski for stimulating discussions regarding conic formulations of nonlinear membranes.

\section*{Data availability}
The datasets generated during and/or analysed during the current study are available in the \texttt{fenics\_optim} repository \cite{fenics_optim_2020} as documented demos.

\bibliographystyle{apalike}
\bibliography{conic_programming_applications.bib}
\end{document}